# Joint Optimization of Charging Infrastructure and Fleet Mix for $CO_2$-Constrained Feeder Services

Haruko Nakao[1], Tai-Yu Ma[2*], Richard D. Connors[3], Francesco Viti[4]

**Haruko Nakao[1]**
[1]University of Luxembourg, Campus Belval, 6, Avenue de la Fonte, Esch-sur-Alzette, Luxembourg
Email: haruko.nakao@uni.lu

**Tai-Yu Ma[2]***
[2]Luxembourg Institute of Socio-Economic Research (LISER), 11 Porte des Sciences, 4366 Esch-sur-Alzette, Luxembourg
Email: tai-yu.ma@liser.lu

**Richard D. Connors[3]**
[3]Luxembourg Institute of Science and Technology (LIST), 5 Avenue des Hauts-Fourneaux, 4362 Esch-sur-Alzette, Luxemburg
Email: richard.connors@list.lu

**Francesco Viti[4]**
[4]University of Luxembourg, Campus Belval, 6, Avenue de la Fonte, Esch-sur-Alzette, Luxembourg
Email: francesco.viti@uni.lu

*Corresponding author

**Abstract**

This study addresses the electrification of demand-responsive feeder services, a form of public transport designed to connect rural and low-demand areas to mass transit hubs. Electrifying demand-responsive transport requires planning the charging infrastructure carefully, considering the trade-offs of charging efficiency and charging infrastructure costs. This study addresses the joint planning of fleet size and charging infrastructure for a demand-responsive feeder service under stochastic demand, given a user-defined $CO_2$ emissions reduction policy. We propose a bi-level optimization model where the upper-level determines charging station configuration given stochastic demand, and the lower-level solves a mix fleet feeder (first and last mile) service routing problem under the $CO_2$ emission and capacitated charging station constraints. An efficient deterministic annealing metaheuristic is proposed to solve the $CO_2$-constrained mixed fleet routing problem. The metaheuristic solves up to 500 requests within 3 minutes, demonstrating the practical applicability of the proposed solution. We applied the model to a real-world case study in Bettembourg, Luxembourg, with two types of electric minibuses and gasoline ones, under different $CO_2$ reduction targets considering rapid (125 kW) and super-fast (220 kW) chargers, given 200 requests per day. The results show that using 24-seat minibuses leads to significant cost savings (-49% on average) compared to that of 10-seat minibuses. Due to their larger battery capacity, charger availability has a smaller impact on the operational costs of 24-seat minibuses. The proposed method provides a flexible tool for joint charging infrastructure and fleet size planning.

## 1. Introduction

Public transport (PT) operates on fixed routes and timetables to meet regional mobility needs. Due to limited resources, PT operators need to navigate between minimising the operational costs and maximising the coverage of service areas. Consequently, the service frequency of fixed-route public transport in low-demand areas is usually insufficient, resulting in low ridership and long waiting times. To address these issues, demand-responsive transport (**DRT**) systems provide a more flexible alternative to allow public transport services to adapt to customers' demand. DRT systems can operate as a door-to-door taxi-like service, a deviation-based service (allowing a part of bus stops optional) (Montenegro et al., 2022), or a feeder service integrating with mass transit systems to reduce its operational costs. Several studies have demonstrated that integrating DRT with mass transit as a feeder service can increase the ridership of public transport and reduce personal car use (Ma et al., 2019; Kumar and Khani, 2021).

In the pursuit of reducing greenhouse gas emissions to 55% below the 1990 level by 2030 (IEA, 2022), electrification of mobility services has become increasingly popular, motivated by $CO_2$ emission reduction targets. At the same time, it is agreed that the modal shift from private cars to mass transit systems is also essential to reach the target level of $CO_2$ reduction (Boehm et al. 2021). We focus on $CO_2$ emissions as it is the primary policy target under EU 'Fit for 55'[1]. A comparison of the economics of batteries, fuel cells, and hydrogen, as well as the production and distribution of different fuels, is beyond the scope of this study. By electrifying Demand Responsive Feeder Services (**DRFS**) (i.e., operating on-demand first/last mile transport services using electric vehicles (**EV**)) would have a huge potential to reduce $CO_2$ emissions directly through the use of EVs and indirectly by increasing the attractiveness (or level of service) of the mass transit. We prioritize electric feeder services because they complement rail, handle multi-passenger trips, and directly replace conventional minibuses. Other modes (bike, e-scooter) are less viable in rural hinterlands and are not viable for all user groups such as the elderly or people with reduced mobility. For the characteristics of shared micromobility systems and their potential to reduce greenhouse gas emissions, car use, and health benefits, the reader is referred to Shaheen and Cohen (2019).

When electrifying their fleet, public transport operators need to carefully consider the trade-offs between investment/operational costs of charging infrastructure and fleet size to serve customer demand. Despite the rapid technological improvement, the battery capacity is not enough for the services to operate continuously for a whole day. Moreover, the charging time of EVs is long, ranging from around 30 minutes to more than 10 hours, depending on the type and power of applied charging technology. In particular, due to the high investment costs of fast charging infrastructure, the number of fast charging stations is still scarce in many cities. This negatively impacts the availability of vehicles due to long charging/waiting times at charging stations. Furthermore, it is difficult to decide on the optimal EV fleet size for DRT systems as vehicle energy consumption rates are easily influenced by environmental factors. Extreme hot weather necessitates the use of the air conditioning system (AC) to ensure passenger comfort, which can lead to an average reduction in driving range by 30-40%, depending on the vehicle size and the AC system's efficiency (Zhang et al, 2018), while the cold weather will reduce battery performance so as driving range (Steinstraeter et al, 2021). Because people's activity locations and travel times fluctuate from day to day (Raux et al., 2016), demand for DRFS exhibits substantial spatiotemporal variability. As a result, even when historic averages are available, the number of trips per vehicle each day and the corresponding energy consumption remain highly uncertain. While seeking for business sustainability, the service operator requires to provide enough fleet to maintain the service level within the acceptable budget, which adds complexity to the planning problem compared to the diesel vehicle operation.

Another factor that influences fleet size optimisation is the configuration of charging infrastructure (i.e., fast chargers in support of overnight slow chargers) (Macrina et al., 2019). If we consider partial or/and fast charging during the operation, the number and location of opportunity charging stations will impact how efficiently the vehicles can include charging activity during the service operation (Yuan et al., 2019 and He et al., 2019). If the number of charging stations is insufficient or the location of charging stations is not optimal, fleet size needs to be larger than necessary to maintain the service level (e.g., waiting time, customer rejection rate). On the other hand, having more than the necessary chargers or charging stations will increase the operational cost as the charging station installation cost is high. Also, once a charging station is installed, it cannot be replaced easily as it is an expensive and time-consuming process. Hence, thoughtful planning and consideration are required to install charging stations in an optimal manner.

[1] https://www.consilium.europa.eu/en/policies/fit-for-55/

Given the high upfront investment cost, the transition from a gasoline-powered fleet to an electric fleet would happen gradually in practice. Some studies started investigating the bus routing scheduling problem with a mix fleet of hybrid vehicles (**HVs**) and Electric Buses (**EBs**) (Rinaldi et al., 2020, and Zhou et al., 2023) or EB replacement plan for transport companies to meet their electrification target in a cost-effective way (Pelletier et al., 2019). While there are many studies in fleet size and charging infrastructure planning for EB or ride-hailing services (Ma and Fang, 2022), there are no flexible modelling approaches and solution algorithms for the fleet size and charging infrastructure planning of DRFS using EVs, taking into account: i) $CO_2$ emission reduction targets, ii) stochastic customer demand, iii) mix fleet. This research issue is particularly relevant as the electrification of the fleet is motivated by $CO_2$ emission reduction, and the level of electrification will be determined by the target $CO_2$ emission reduction and the investment cost to achieve the target. We assume plug-in recharging from the grid. Section 4 provides specifications for the scenarios covered. Note that we do not consider battery swapping. There are some very limited examples of battery swapping deployed in Europe for passenger EVs, but to our knowledge, battery swapping is not yet used in commercial electric bus or minibus operations, due to high mechanical complexity, large battery pack size, and lack of standardization across (mini)bus manufacturers.

With focus on low-demand rural contexts, this study aims to develop a generalized model for transit feeder service operation to plan the electrification of their fleet considering the environmental impact, especially focusing on the trade-off between $CO_2$ emissions and the operational and infrastructure installation cost associated with the electrification. In Europe the dominant type of electric buses are battery-electric buses (BEBs) ), whereas the PHEV technology is gradually being phased out. The ICCT report[2] shows European bus sales were 36% battery electric and only 1% PHEV. In Europe. We therefore limit ourselves to considering the mixed fleet size planning problem with two powertrain technologies: BEBs and gasoline buses (GBs). The developed methods would support the practitioners in their strategic and tactical decisions during and beyond the transition period. The contributions of this paper are as follows.

a) We develop new integrated planning approach to jointly optimize the fleet size and charging infrastructure of DRFS service. The approach provides a flexible modelling framework allowing the operator to configure mixed fleet size and $CO_2$ emission reduction target. A bi-level optimisation model (two-stage stochastic program) is developed where the upper-level problem configures the location, number and type of chargers, and the lower-level problem determines the best mix fleet composition to meet customer demand and the $CO_2$ reduction target. An efficient iterative solution algorithm is proposed to solve the bi-level problem efficiently.
b) A mixed integer linear programming (**MILP**) model is proposed for mixed bus fleet routing optimization of the feeder service with $CO_2$ emission and vehicle routing constraints. The charging scheduling sub-problem considers the charging station capacity limits with a partial recharge policy and different charger types. A series of valid inequalities is proposed to tighten the MILP formulation and solve the problem more efficiently.
c) We develop an efficient deterministic annealing (DA) metaheuristic to solve the $CO_2$-emission-constrained mixed fleet feeder service routing problem efficiently. A set of tailored operators are developed to address complicated mixed fleet utilisation with $CO_2$ emission constraints and charging schedule synchronisation to meet charging station capacity. The results show that the proposed algorithm can obtain near-optimal solutions fast (e.g., around 3 minutes or less for 500 requests with 50% or 90% $CO_2$ reduction targets).
d) The application of the developed model is conducted for a real-world case study in the Luxembourg Bettembourg area. Different scenarios were designed based on two types of EVs (10-seat and 24-seat EVs), fast and superchargers, as well as 50% and 90% $CO_2$ reduction targets. The results show that with 200 daily requests, using 24-seat EVs would lead to significant system cost savings (-49% on average) due to a smaller required fleet size and charging infrastructure. The results show that the proposed model and algorithm can obtain an optimised (mixed) fleet size and charging station configuration under the desired $CO_2$ reduction target under stochastic demand.

[2] https://theicct.org/publication/r2z-eu-hdv-market-development-quarterly-jan-sept-2024-dec24/

The remainder of this paper is organized as follows. In Section 2, the relevant literature is discussed. Section 3 provides the problem statement, model formulation, and solution algorithm. Section 4 explains the test instance generation and presents the procedure and computational results for the parameter tuning of the DA metaheuristic. Then, we compare the performance of the DA algorithm with the exact solutions obtained from a commercial MILP solver using a series of test instances. To evaluate the efficiency of the DA algorithm, we test it on a series of larger test instances with up to 500 requests, given different $CO_2$ emission targets. Section 5 is dedicated to a real-world case study. First, the description of the case study and considered scenarios are presented. It follows the analysis of the computational results and discussion on the managerial insights. Finally, conclusions are drawn, and future extensions are discussed in Section 5.

## 2. Literature review

There are extensive studies focusing on fleet size and charging infrastructure planning for fixed-route bus systems using EBs only or with a mixed fleet of EBs and GBs. For the DRT system, Ma and Fang (2022) provide a systematic review of the models and algorithms for strategic, tactical, and operational decisions. In this section, we review the relevant literature related to fleet size and charging infrastructure planning for fixed-route bus and DRT systems.

Electrifying **fixed-route bus systems** relies on sufficient charging infrastructure capacity to meet the charging needs of the fleet. As the charging infrastructure configuration depends on the charging demand, which is affected by the fleet size and the vehicle's characteristics (e.g., energy consumption rate, battery capacity), existing studies jointly optimise the location/configuration of charging stations and the composition of EBs and charging schedules. Rogge et al. (2018) proposed a MILP for the strategic planning of EB fleet composition (two types of EBs) and the number of homogeneous chargers at the single depot to serve a set of scheduled trips. The proposed genetic algorithm is applied for a case study involving a couple of bus routes in Aachen, Germany. Alamatsaz et al. (2024) propose both arc-based MILP and path-based integer linear programming approaches to jointly optimize the fast charging infrastructure and charging schedules of EBs. Their model assumes full recharging and a fixed charging time for each charging event. The authors developed an exact branch-and-price (B&P) algorithm to solve the model on a set of test instances with up to 70 bus trips and 8 candidate charging station locations. Nath et al. (2024) developed a comprehensive MILP model to optimize charging station location, trip assignment and charging scheduling. The model considers several practical aspects, including grid capacity constraints, time-of-use energy prices, and partial recharging. An iterated local search algorithm is developed to solve the integrated planning problem using real-world GTFS datasets with multiple depots, hundreds of buses, and thousands of bus trips and stops.

Other studies consider joint fleet composition and charging infrastructure planning with EBs and GBs. For example, Hu et al. (2022) optimised the location of fast chargers to be installed at selected bus stops and the battery capacity of EBs with a given fleet size as well as provided the bus charging schedules. They minimised the overall system costs also considering the time-varying electricity price as user's disutility, such as extra waiting time at the stop. The advantage of this joint optimisation approach is a more cost-efficient investment. However, this might need to solve a more complicated optimisation problem for which efficient solution algorithms remain a research challenge. Lee et al. (2021) developed a MILP for EB fleet size and charging infrastructure planning to configure fleet size, battery capacity and charging infrastructure for the electrification of a single bus line. A scenario-based approach is applied to model the impact of uncertain bus energy consumption. A two-stage sequential algorithm (determining feasible configuration of the fleet size and number of chargers at the first stage, then minimising the charging cost at the second stage) is developed and applied for a case study in Seogwipo city, Korea, using bus operation data. An (2020) proposed a stochastic integer program for optimising fleet size and charging station location with stochastic charging demand and time-of-use electricity price. The stochastic integer program is approximated as a MILP by the sampling approach under different demand scenarios. A Lagrangian relaxation approach is developed to solve the MILP problem. The methodology is applied for a case study for a bus network case in Melbourne, Australia.

Recent studies address both charging station locations and its configurations as well as fleet compositions and schedules (Guschinsky et al., 2021; Battaïa et al., 2023). Guschinsky et al. (2021) proposed mixed integer programming models to optimise the location of charging stations, number of transformers connected to charging stations, number of chargers per station, and number of connectors to buses per charger as well as the fleet composition (EB and GB) and its assignment to fixed bus routes. Yildrim and Yildiz (2021) investigated the impact of different types of charging technology (i.e., fast charger and dynamic wireless power transfer and its configuration (i.e., coverage area of DWPT, spatial distribution of fast charger among multiple depots and their charging power) to the total electrification cost assuming the heterogeneous EB fleet. A heuristic based on the column generation approach is proposed and applied for a real-world bus networks case study in Istanbul. Soltanpour et al. (2023) developed a joint fleet size and charging infrastructure planning model for a mixed fleet. The model assumes a full recharging policy and time-dependent energy prices. To ensure compliance with scheduled timetables during bus charging operations, the optimal number of backup buses is determined. Scenario-based approaches are used to analyze the effect of charging powers (150 kW, 200 kW, and 350 kW) and bus size (20 seats and above) on total costs. Sen et al. (2024) developed a MILP model that jointly optimizes charging schedules and bus trip assignments given fixed charging infrastructure and a mixed fleet of EBs and GBs while considering dynamic electricity prices. The model uses transit blocks (sets of sequential trips assigned to a bus during the day) as analytical units. Based on these units, the charging scheduling model is developed to minimize total operating costs.

For DRT systems, recent studies have embedded an eco-oriented routing strategy into the routing strategy to optimize its energy efficiency. For instance, Li et al. (2022) proposed an optimization framework for eco-oriented DRT system routing design. They propose a MILP model that integrates fuel consumption and $CO_2$ emission constraints to optimize vehicle routing. Then, a dynamic module adapts the MILP solution in response to road conditions. The MILP model is formulated for day-ahead planning under the assumption that all ride requests and road traffic conditions are known. Shi et al. (2025) developed a bi-objective MILP model for the mixed-fleet vehicle routing problem in collaborative logistics, which is a joint operation between logistics operators who share delivery tasks and transportation resources to improve efficiency and sustainability. The first objective is to maximize total profit by reducing transport costs and charging expenses. The second objective is to minimize $CO_2$ emissions from internal combustion engine vehicles. $CO_2$ emissions from vehicles are proportional to vehicle load weight and travel distance. The capacity of charging stations is assumed to be unlimited. Deka et al. (2023) provide a detailed review of progress in the sustainability of DRT services.

The joint fleet size and charging infrastructure optimisation of DRT systems contains the same nature of the trade-off between fleet size and charging infrastructures as the fixed-route bus systems. However, as the routes and schedules of on-demand services are not fixed, the impact of changes in one value (e.g., fleet size or location of charging stations) could completely renew the vehicle's routes and schedules. Therefore, it is more computationally expensive to find the optimal fleet size and charging infrastructures compared to fixed-route bus systems. As a result, many studies often focus on either the optimal fleet sizes for EV operated on-demand services (e.g., Hiermann et al., 2016) or charging configurations (Ma and Xie, 2021; Li et al., 2022; Alam and Guo, 2022). Recently, Zhang et al. (2025) proposed a joint optimization approach for fleet operations and charging infrastructure planning of ride-hailing systems using autonomous vehicles.

There are few studies that jointly optimised the fleet size and charging infrastructure of on-demand services simultaneously. Some studies take an agent-based simulation approach to estimate the required fleet size and evaluate most economic charging infrastructure configurations for shared electric autonomous vehicles (e.g., Chen et al., 2016). Yang et al. (2023) proposed the simulation framework for electric autonomous mobility-on-demand to design the fleet size and charging facility configuration by integrating the mathematical model to 1) assign a vehicle to request, 2) relocate empty vehicles, and 3) assign a vehicle to the charging facility. They also considered endogenous congestion. However, the decision problem is to determine the fleet size of homogeneous EVs instead of mixed fleets, the latter being more complicated to configurate. Other studies model the problem as a MILP to find the optimal solutions. Zhang et al. (2020) is one of the first to develop a MILP for modelling joint fleet sizing and charging infrastructure optimisation problem. They assume that autonomous electric vehicles are operated to serve passenger and goods transportation. All the nodes in the network were assumed to be potential charging stations, and they optimised the number of chargers per station. These assumptions limit the application of this model to intercity trips where the node and link represent cities and corridors between them. Paudel and Das (2022) developed a profit-maximizing MILP model to determine the fleet size, the location of the charging hub, and the number of chargers per hub for ride-hailing services operated by shared autonomous electric vehicles. The objective is to maximize the operator's profit, but not all demand is necessarily served in this model. A robust model is formulated to deal with the uncertainty of transportation demand. The authors generated multiple demand scenarios (e.g., 100) given the mean and standard deviation. Then, the fleet and charging infrastructure configurations are determined given a targeted level of demand coverage (e.g., all 100 scenarios will be covered or up to mean demand will be covered).

Except for Paudel and Das (2022), existing studies do not consider the stochasticity of demand, which is important to provide robust solutions. Besides, all three papers cited above have assumed a fully electric fleet, which does not allow the consideration of the mixed fleet as well as a $CO_2$ emission reduction target. Hence, these approaches fail to reveal to what extent charging infrastructure configuration impacts upon service performance under different $CO_2$ reduction targets. This study aims to address these research gaps, providing a new joint fleet size and charging infrastructure optimization method under stochastic customer demand given the configurable $CO_2$ emission reduction target, charging station capacity constraints, and partial recharge.

Appendix A1 summarises the characteristics of the related studies cited in this literature review, and highlights the novelty of this study.

## 3. Methodology

### *Notations*

| | |
|---|---|
| ***Sets and variables*** | |
| $0, N+1$ | Two instances of the depot |
| $P$ | Set of pick-up vertices (locations) associated with demand scenario $\xi$, $P = \{1, \dots, n\}$, ($\xi$ is dropped) |
| $D$ | Set of drop-off vertices (locations) associated with demand scenario $\xi$, $D = \{n+1, \dots, 2n\}$ ($\xi$ is dropped) |
| $P_{N+1}$ | $P_{N+1} = P \cup \{N+1\}$ |
| $D_0$ | $D_0 = D \cup \{0\}$ |
| $T$ | Set of dummy vertices associated with transit stations, $T \subset P \cup D$ |
| $S$ | Set of physical chargers, i.e. $S = \{1, \dots, \bar{S}\}$ |
| $S'$ | Set of dummy charger vertices, $S' = \cup_{o \in S} S'_o$, where $S'_o$ is the set of dummy charger vertices for a physical charger $o, o \in S$. |
| $A$ | Set of arcs |
| $\mathcal{M}$ | Set of vehicle types (electric and gasoline) with $\mathcal{M} = \{e, c\}$ |
| $K$ | Set of vehicles, i.e. $K = K^e \cup K^c$ where $K^e$ is the set of electric vehicles and $K^c$ is the set of conventional vehicles with the combustion engine. We assume that the conventional gasoline/electric vehicles are the same type. The set of vehicles can be defined as $K = \{K^m\}_{m \in \mathcal{M}}$, where $K^m = \{1, \dots, n_m\}$, assuming $n_m$ is enough to serve all requests. |
| $W$ | Set of candidate locations to install charging stations |
| $H$ | Set of charger type |
| $B_i^k$ | Beginning time of service of vehicle $k$ at vertex $i$ |
| $V$ | Set of vertices except for the depot, i.e. $V = P \cup D \cup S'$ |
| $V_0, V_{N+1}, V_{0,N+1}$ | $V_0 = V \cup \{0\}, V_{N+1} = V \cup \{N+1\}, V_{0,N+1} = V \cup \{0, N+1\}$ |
| $Q_i^k$ | Load of vehicle $k$ at vertex $i$, $\forall k \in K$ |
| $q_i$ | Passenger load at vertex $i$, $\forall\, i \in P$ |
| $E_i^k$ | State of charge of vehicle $k$ at vertex $i$, $\forall k \in K^e$ |
| | |
| ***Parameters*** | |
| $c_{ij}$ | Distance from vertex $i$ to vertex $j$ |
| $t_{ij}$ | Vehicle travel time from vertex $i$ to vertex $j$ |
| $L_i$ | Maximum ride time of request $i$ being the direct ride time multiplied by a detour factor (e.g., 1.5) |
| $u_i$ | Service time at vertex $i$ |
| $Q_{max}^k$ | Capacity of vehicle $k$ |
| $[e^i, l^i]$ | Earliest and latest starting times of service associated at vertex $i$ |
| $\Gamma(\xi)$ | Total $CO_2$ emission for serving demand scenario ξ using GVs only |
| $\pi$ | Targeted user-defined $CO_2$ emission reduction rate $0 \leq \pi \leq 1$ |
| $E_{min}^k, E_{max}^k, E_{init}^k$ | Minimum, maximum, initial states of charge of electric vehicle $k, \forall k \in K^e$ |
| $\alpha_s$ | Charging rate of charger $s \in S'$ |
| $\beta^k$ | Energy consumption rate for vehicle k, depending on the type of vehicles (kWh/km). For conventional vehicles, we assume an average value of $\beta^c$(l/km). |
| $f_0$ | Charging station investment cost, converted to euros/day |
| $f_h$ | Unitary cost of type $h$ (fast/rapid) charger, converted to euros/day |
| $\bar{f}_k$ | Unitary cost of slow charger, converted to euros/day. $\bar{f}_k = 0$ if vehicle $k$ is conventional vehicle |
| $f_k$ | Unitary purchase cost of vehicle $k$, converted to euros/day |
| $\vartheta_h$ | Charging power per hour for the type-$h$ charger (kW) |
| $\omega_w$ | The power supply constraints at the charging station $w$ |
| $\lambda^k$ | Energy price for vehicle $k$, $\forall \lambda^k \in \{\lambda^c, \lambda^e\}$, where $\lambda^c$ (euro/l) and $\lambda^e$(euro/kWh) is the energy price for conventional and electric vehicles, respectively |
| $\xi$ | Demand scenario for the lower-level problem (vehicle routing of feeder service), $\forall \xi \in \mathcal{E}$ |
| $b$ | Total number of types of electric vehicles |
| $\theta^k$ | $CO_2$ emission rate per energy consumption of vehicle k; gasoline vehicles: $\theta^c$ (kg/l), and electric vehicle: $\theta^e$ (kg/kWh), i.e., $\theta^k \in \{\theta^c, \theta^e\}$ |
| $p_\xi$ | Probability of the demand scenario $\xi$ with $\sum_{\xi \in \mathcal{E}} p_\xi = 1$ |
| $M_0, M_1, M_2, M_3$ | Sufficiently large positive numbers |
| ***Decision variable*** | |
| $\tilde{y}_w$ | 1 if a charging station is built in the charging location candidate $w$, and 0 otherwise |
| $y_{hw}$ | The number of type-$h$ chargers to be installed at the charging station $w \in W$ |
| $x_{ij}^k$ | 1 if arc $(i, j)$ is travelled by vehicle $k$, and 0 otherwise |
| $\tau_s^k$ | Charging duration for vehicle $k$ at charger $s$, $s \in S'$ |

3.1. Optimisation model for the joint fleet size and charging infrastructure planning

We consider the joint fleet size and charging infrastructure planning problem for a DRFS operated by a Transport Network Company (**TNC**) in a rural area. The DRFS operates a fleet of minibuses on a set of pre-defined virtual (optional) bus stops that densely cover the service area to ensure all potential customers' origins/destinations are within reasonable walking distance. Customers submit their trip requests specifying the pick-up/drop-off locations and arrival times. Destinations are either bus stops or a train station where scheduled transit service is offered. Assuming the TNC owns enough GVs to cover the current demand, they are interested in the electrification of (part of) the fleet to meet a $CO_2$ emission reduction target. **Throughout the remainder of this paper, the EVs being considered are electric minibuses**. The $CO_2$ emission reduction target can be determined by policy (e.g., the EU's goal of reducing GHG emissions by 55% by 2030 compared to 1990 levels, according to the International Energy Agency (IEA, 2022)), stakeholder preference (e.g., Luxembourg City plans to phase out combustion engine buses by 2030[3]), or the operator's commitment to sustainability. However, some trade-off decisions arise from the cost of fleet electrification and the variability of daily customer demand. We consider total electrification costs to include charging infrastructure investment costs, fleet acquisition costs of EVs, and energy consumption costs to serve the demand, all converted into per diem costs assuming a vehicle lifetime (e.g., 10 years). Customer demand is assumed stochastic following a uniform probability distribution $p$. In practice, an empirical probability distribution can be used instead if historical ride data is available. Note that we only model typical weekday demand and exclude weekend demand patterns. The primary goal of our strategic planning model is to ensure peak weekday demand and the necessary fleet size and charging infrastructure. Given a target $CO_2$ reduction rate, the TNC aims to minimize the electrification cost of DRFS while satisfying customer demand within some service levels (ensuring all customer requests when demand is normal). We assume that there are several candidate electric vehicle types and candidate locations for building charging stations with fast chargers, which allow EVs to be recharged during their operation and return to service quickly. One charging station can have multiple chargers, which allows multiple EVs to be charged simultaneously. The number of chargers per charging station is limited so that the maximum total power usage from the grid does not exceed its capacity. To maximize the mixed fleet DRFS's performance, the charging infrastructure configuration decisions (i.e., location, number of chargers, and their types) are decided jointly with the EV acquisition (i.e. decide the number of EVs to be purchased for each vehicle type), while meeting constraints of the local power grid. The numerical example in Section 4 provides examples of relevant EV and charger specifications.

[3] https://delano.lu/article/35m-to-electrify-77-regional-b

The problem is formulated as a bi-level optimisation problem where the overall electrification system costs of the feeder service is to be minimised. The upper level determines the number of charging stations and chargers installed at each charging station while the lower level determines the fleet size and composition of DRFS to meet scenario-based customers' demand with targeted $CO_2$ emission constraint. The upper-level problem is formulated as follows.

(***L1***)

$$min\, Z_U = \sum_{h\epsilon H}\sum_{w\in W}(f_h y_{hw} + f_0\tilde{y}_w) + \bar{Z}_L(\boldsymbol{y},\widetilde{\boldsymbol{y}}) \tag{1}$$

subject to:

$$\sum_{h\in H} y_{hw} \le M_0\tilde{y}_w, \qquad \forall w \in W \tag{2}$$

$$\sum_{h\in H} y_{hw}\vartheta_h \le \omega_w, \quad \forall w \in W \tag{3}$$

$$\tilde{y}_w \in \{0,1\}, \qquad \forall w \in W \tag{4}$$

$$y_{hw} \in Z_0^+, \quad \forall h \in H, w \in W \tag{5}$$

The objective function (1) minimises the total system costs (including charging infrastructure investment, fleet acquisition, and energy consumption of service trips) regarding fleet electrification. The first term indicates the total installation and acquisition costs of charging stations ($\tilde{y}_w$) and a set of fast/rapid chargers ($y_{hw}$). The second term $\bar{Z}_L(\boldsymbol{y},\widetilde{\boldsymbol{y}})$ is the total expected costs of the fleet use costs (converted daily vehicle use costs considering acquisition and maintenance) and energy consumption costs to serve customers' demand. Equation (2) ensures that chargers can be installed at a location $w$ when it is open. Equation (3) ensures that the number of installed chargers at a location $w$ cannot exceed its electric power supply constraints. We assume that the type of chargers installed on the same site is homogeneous. Eqs. (4)-(5) are domain variables.

The term $\bar{Z}_L(\boldsymbol{y},\widetilde{\boldsymbol{y}})$ is obtained by solving the following lower-level stochastic optimization problem to minimize expected fleet size and operational costs given a set of customers' demand scenarios.

(***L2***)

$$\bar{Z}_L(\boldsymbol{y},\widetilde{\boldsymbol{y}}) = \frac{1}{|\mathcal{E}|}\sum_{\xi\in\mathcal{E}} p_\xi Z_L(\xi) \tag{6}$$

St.

$$S \in \Omega(\boldsymbol{y},\widetilde{\boldsymbol{y}}) \tag{7}$$

We assume historical customer demand data is available and is used for estimating the probability distribution $p_\xi$ of stochastic demand scenarios $\xi \in \mathcal{E}$. As a result, the expected cost over $\mathcal{E}$ can be estimated by Eq. (6), where $Z_L(\xi)$ is the objective function value of the Feeder (first and last mile) Service Mix Fleet

Routing Problem (**FS-MFRP**), given a demand scenario $\xi$. This is formulated as a MILP problem of Eqs. (8)-(36). To get the expected $Z_L(\xi)$, the FS-MFRP needs to be solved $|\mathcal{E}|$ times. $\Omega(\boldsymbol{y},\widetilde{\boldsymbol{y}})$ in Eq. (7) is a feasible charging facility configuration for FS-MFRP. The objective of FS-MFRP (Eq. (8)) takes into account the costs of the complementary charging connectors of (slow) chargers to be installed at the depot (assuming each vehicle has their respective overnight slow charger(connector) at the depot) for ensuring overnight charging availability for each vehicle, given a charging station configuration in $(\boldsymbol{y},\widetilde{\boldsymbol{y}})$. Note that most DRT operators (e.g., Uber and Lyft, among many others) use dedicated applications to manage and collect users' ride requests and information. Since our model is used for strategic planning, it is reasonable to assume average or smoothed demand and neglect extreme demand peaks (outliers). When demand data is incomplete or unavailable, operators can use other data sources, such as mobility survey data or synthetic populations and their mobility patterns, to estimate customer demand. Our approach assumes that historical ride data is available in order to obtain the empirical probability distribution of demand scenarios. This scenario-based sampling approach has been applied to a similar joint fleet size and charging infrastructure planning optimization problem for electric buses (An, 2020). Another research avenue for the future extension of this study is using a robust optimization approach (Paudel and Das, 2022). We acknowledge that data privacy issues must be addressed by operators when collecting users' ride data and information. This implies strict compliance with GDPR regulations, anonymization of collected data, and implementation of a data management plan for data protection.

Given a charging infrastructure configuration, the FS-MFRP problem is formulated as follows, assuming customer demand $\xi \in \mathcal{E}$ to be met by a mixed fleet of vehicles, given a user-defined $CO_2$ emission reduction target $\pi$ with $0 < \pi < 1$ with respect to the total emission when using a fleet of GVs.

(***FS-MFRP***)

$$Z_L(\xi) = min \sum_{k\in K}\sum_{(i,j)\in A} \lambda^k \beta^k c_{ij}\, x_{ij}^k + \sum_{k\in K}\sum_{j\in V} (f_k + \bar{f}_k)\, x_{0j}^k \tag{8}$$

Subject to:

$$\sum_{k\in K^c}\sum_{(i,j)\in A} \theta^k \beta^k c_{ij} x_{ij}^k \leq (1-\pi)\Gamma(\xi) \tag{9}$$

$$\sum_{j\in P_{N+1}\cup S'} x_{0j}^k = \sum_{i\in D_0\cup S'} x_{i,N+1}^k = 1, \qquad \forall k \in K \tag{10}$$

$$\sum_{k\in K}\sum_{j\in V_0} x_{ij}^k = 1, \qquad \forall i \in P \tag{11}$$

$$\sum_{j\in P_{N+1}} x_{sj}^k \leq 1, \qquad \forall k \in K^e, s \in S' \tag{12}$$

$$\sum_{j\in V_{0,N+1}} x_{sj}^k = 0, \qquad \forall k \in K^c, s \in S' \tag{13}$$

$$\sum_{j\in V} x_{ij}^k - \sum_{j\in V_{N+1}} x_{n+i,j}^k = 0, \qquad \forall k \in K, i \in P \tag{14}$$

$$\sum_{j\in V_0} x_{ji}^k - \sum_{j\in V_{N+1}} x_{ij}^k = 0, \qquad \forall k \in K, i \in P \cup D \tag{15}$$

$$Q_j^k \geq Q_i^k + q_j - M_1\left(1 - x_{ij}^k\right), \qquad \forall k \in K, i,j \in V_{0,N+1} \tag{16}$$

$$Q_j^k \leq Q_i^k + q_j + M_1\left(1 - x_{ij}^k\right), \qquad \forall k \in K, i,j \in V_{0,N+1} \tag{17}$$

$$0 \leq Q_i^k \leq Q_{max}^k, \qquad \forall k \in K, i \in V_{0,N+1} \tag{18}$$

$$B_j^k \geq B_i^k + u_i + t_{ij} - M_2(1 - x_{ij}^k), \qquad \forall k \in K, i \in V_{0,N+1}\backslash S', j \in V_{0,N+1} \tag{19}$$

$$B_j^k \geq B_s^k + u_s + \tau_s^k + t_{sj} - M_2(1 - x_{sj}^k), \qquad \forall k \in K^e, s \in S', j \in P_{N+1} \tag{20}$$

$$e^i \leq B_i^k \leq l^i, \qquad \forall k \in K, i \in P \cup D \tag{21}$$

$$L_i \geq B_{n+i}^k - \left(B_i^k + u_i\right), \qquad \forall i \in P, k \in K \tag{22}$$

$$E_0^k = E_{init}^k, \qquad \forall k \in K^e \tag{23}$$

$$E_{min}^k \leq E_i^k \leq E_{max}^k, \qquad \forall k \in K^e, i \in V \tag{24}$$

$$E_j^k \geq E_i^k - \beta^k c_{ij} - M_3\left(1 - x_{ij}^k\right), \qquad \forall k \in K^e, i \in V_0\backslash S', j \in V_{N+1} \tag{25}$$

$$E_j^k \leq E_i^k - \beta^k c_{ij} + M_3\left(1 - x_{ij}^k\right), \qquad \forall k \in K^e, i \in V_0 \backslash S', j \in V_{N+1} \tag{26}$$
$$E_j^k \geq E_s^k + \alpha_s \tau_s^k - \beta^k c_{sj} - M_3\left(1 - x_{sj}^k\right), \qquad \forall k \in K^e, s \in S', j \in P_{N+1} \tag{27}$$
$$E_j^k \leq E_s^k + \alpha_s \tau_s^k - \beta^k c_{sj} + M_3\left(1 - x_{sj}^k\right), \qquad \forall k \in K^e, s \in S', j \in P_{N+1} \tag{28}$$
$$v_s = \sum_{j \in P_{N+1}} \sum_{k \in K^e} x_{sj}^k, s \in S' \tag{29}$$
$$v_h \leq v_l, \forall h, l \in S_o', o \in S, h < l \tag{30}$$
$$\sum_{k \in K^e} B_h^k \geq \sum_{k \in K^e} B_l^k + \sum_{k \in K^e} \tau_l^k - M_2(2 - v_h - v_l), \forall h, l \in S_o', o \in S, h < l \tag{31}$$
$$\tau_s^k + B_s^k \leq M_2 \sum_{j \in P_{N+1}} x_{sj}^k, \forall s \in S', k \in K^e \tag{32}$$
$$v_s \leq 1, s \in S' \tag{33}$$
$$x_{ij}^k \in \{0,1\}, \qquad \forall k \in K, i, j \in V_{0,N+1} \tag{34}$$
$$B_i^k \geq 0, \forall k \in K, i \in V_{0,N+1} \tag{35}$$
$$\tau_s^k \geq 0, \forall k \in K^e, s \in S' \tag{36}$$

The FS-MFRP is modelled on a complete graph $G(V, A)$ where $V$ is a set of vertices, and $A$ is a set of arcs. Each arc $(i, j)$ is associated with a travel distance $c_{ij}$ and travel time $t_{ij}$. For simplification, the energy consumption on arcs is a linear function of the travel distance and the average energy consumption rate of vehicles. Different from the classical static Dial-a-Ride problem (**DARP**), we consider typical weekday customers' demand patterns following some probability distribution $\mathbb{P}$. We denote $\xi$ a (demand) scenario as a sampling of demands from $\mathbb{P}$ for a planning horizon $T$ (i.e. one day). For each day, a set of $n_\xi$ requests are generated. Each request is characterised by a number of customers with the same pickup and drop-off locations and the same associated time windows. The time window defines the earliest and latest time interval that a bus service can start at the location of a pickup/drop-off/transit station (Cordeau and Laporte, 2003). Time windows are set based on transit timetables and the buffer time before and after train departures/arrivals to ensure reasonable transfer/waiting times at transit stations. It should be noted that each request may contain one or more passengers heading for the same train station with the same desired arrival time, generated based on known probability distributions from the operator. For the outbound requests (first-mile requests), the latest bus arrival time at the transit station node $i \in T$ is set as $l^i = t^{dep(i)}$ where $t^{dep(i)}$ is the scheduled departure time of the train. The earliest bus arrival time is $e^i = l^i - \delta$, where $\delta$ is the transfer buffer time (e.g. 10-15 minutes). For the inbound requests (last-mile requests), the time windows are associated with their pickup locations (transit stations). The earliest arrival time is $e^i = t^{arr(i)}$ where $t^{arr(i)}$ is the scheduled arrival time of the train. The latest bus arrival time is defined by extending $e^i$ with a transfer buffer time $\delta$, i.e., $l^i = e^i + \delta$.

The widths of time windows are system parameters to be specified by the operator. For example, if users indicate their desired departure time of trains at a transit station is 9:00, the corresponding latest arrival time windows can be set as a couple of minutes (buffer time) before 9:00 to allow users to access their train. Similarly, when users get off a train at time $t$, the earliest and latest arrival times can be set up around t. We assume a hard time window constraint so as to commit to reliable feeder service. In practice, soft time windows can also be applied to penalise early and/or late arrivals. Based on this setting, each request is associated with one pickup and one drop-off vertices. In case two different requests have the same physical pickup/drop-off locations, (duplicate) dummy vertices are generated to allow each vertex to be visited once for modelling the problem.

The objective function (8) minimizes the overall energy consumption costs, fleet acquisition/maintenance costs $(f_k)$, and the costs of slow overnight chargers $(\bar{f}_k)$ -assuming one slow charger (charging plug) is installed for each electric vehicle. Equation (9) is the targeted $CO_2$ emission constraint. The RHS (right-hand side) of Eq. (9) is the targeted $CO_2$ emission, where $\Gamma(\xi)$ is the reference $CO_2$ emission level when serving the demand scenario $\xi$ using exclusively GVs. Given a pre-defined percentage reduction rate $(\pi)$, the targeted total $CO_2$ emission amount varies from one (demand) scenario to another. As a result, we need to solve the FS-MFRP twice for each scenario. First, we solve FS-MFRP using GVs only without Eq. (9) to obtain the reference emission level $\Gamma(\xi)$. Then, we introduce Eq. (9) and solve the problem anew by considering a mixed fleet of vehicles and the reduction rate $\pi$. Constraint (10) ensures that each vehicle starts and ends its service at the depot. Constraint (11) ensures that each pickup location is visited once. Constraint (12) states that a vehicle can be recharged at most once at each dummy charger node, while constraint (13) forbids a gasoline-powered vehicle from visiting the electric charging station. Constraint (14) ensures that a customer is picked up and dropped off by the same vehicle. The flow conservation at each node is guaranteed by constraint (15). Constraints (16)-(18) state the changes of (customer) load at vehicle stops and transit stations as well as the load capacity constraints. The beginning time of service at pickup and drop-off nodes and charger nodes is ensured by Eqs. (19) and (20), respectively. Constraint (21) expresses the time window constraints associated with transit stations. Eq. (22) states that a customer's ride time cannot exceed a maximum value, parameterised by a detour factor with respect to the direct ride time between the passenger's pickup and drop-off locations. Constraints (23) to (28) state that vehicles' state-of-charge (SOC) is bounded and the changes of SOC when traversing an arc or recharging at charging stations. Constraints (29) to (32) ensure that two vehicles cannot use the same charger simultaneously, and one must wait until the other finishes its charging. Constraints (33) to (35) define the domain of all decision/auxiliary variables. Note that the model can be extended to the cases with multiple depots and multiple types of electric vehicles without difficulty (Braekers et al., 2014).

### 3.2. FS-MFRP problem preprocessing and valid inequalities

To tighten the model, we adapt the following preprocessing rules of time-window tightening and arc elimination (Cordeau, 2006). Let $n$ denote the number of requests. For each pair of requests $(i,j), \forall i,j \in P$, we check whether they are compatible by verifying the feasibility of the following paths (Cordeau, 2006): $\{i,j,n+i,n+j\}, \{i,j,n+j,n+i\}, \{j,i,n+i,n+j\}, \{j,i,n+j,n+i\}, \{i,n+i,j,n+j\}, \{j,n+j,i,n+i\}$. If two requests $(i,j)$ are incompatible (i.e. they cannot be served by the same vehicle), the arcs $\{i,n+i\}$ and $\{j,n+j\}$ are eliminated. The model is further tightened by eliminating the symmetry of used vehicles for each vehicle type. We add the following strong valid inequalities.

$$\sum_{j\in V|(0,j)\in A} x_{0j}^{k} - \sum_{j\in V|(0,j)\in A} x_{0j}^{k+1} \geq 0, \qquad \forall k \in K^m \backslash n_m, \forall m \in \mathcal{M} \tag{37}$$

$$M_0 \sum_{j\in V|(0,j)\in A} x_{0j}^{k} \geq \sum_{i,j\in V|(i,j)\in A} x_{ij}^{k}, \qquad \forall k \in K^m, \forall m \in \mathcal{M} \tag{38}$$

Equation (37) ensures that for each type of vehicle $m$, the $k$-th vehicle of that type is used before the *k+1*-th vehicle of the same type. $n_m$ is the number of vehicles of type $m$. Equation (38) states that if a vehicle does not leave the depot, this vehicle is not used to serve the requests. Note $M_0$ is set as the number of arcs in $\{(i,j) \in A | i,j \in V\}$, and $A$ is the set of arcs after trimming off infeasible arcs based on the above rules.

For charging operations, we adapt the valid inequalities when visiting a charger violates the time window constraints. Following Bongiovanni et al. (2019), if both paths $\{n+i,s,j\}$ and $\{n+j,s,i\}$ with $i,j \in P, s \in S'$ are infeasible, the following strong valid inequality is added.

$$\sum_{k\in K^e} \left(x_{n+i,s}^{k} + x_{sj}^{k} + x_{n+j,s}^{k} + x_{si}^{k}\right) \leq 1 \tag{39}$$

### 3.3. Solution algorithm

This bi-level optimization problem is complicated to solve (requiring an iterative procedure to evaluate the total system costs of feasible charging infrastructure configuration), as the lower-level problem is a variant of the mixed fleet DARP (an NP-hard problem). The solution framework is shown in Figure 1. For each charging infrastructure configuration, the lower-level problem is solved for multiple demand scenarios (sample average approximation) from which the expected DRFS electrification cost is estimated (lines 5-8 in Algorithm 1). Then, the minimum expected electrification cost among different charging configurations, given the number of chargers, is compared with the one with the greater number of chargers. If the minimum expected electrification cost decreases by increasing the number of chargers, the simulation is extended to the greater number of chargers. Otherwise, the simulation is terminated with the current number of chargers (lines 14-15 in Algorithm 1). The electrification cost is further compared among different allocations of chargers at the charging stations. The solution with the lowest total electrification cost among them is determined as the optimal solution. The solution algorithm is shown in **Algorithm 1**.

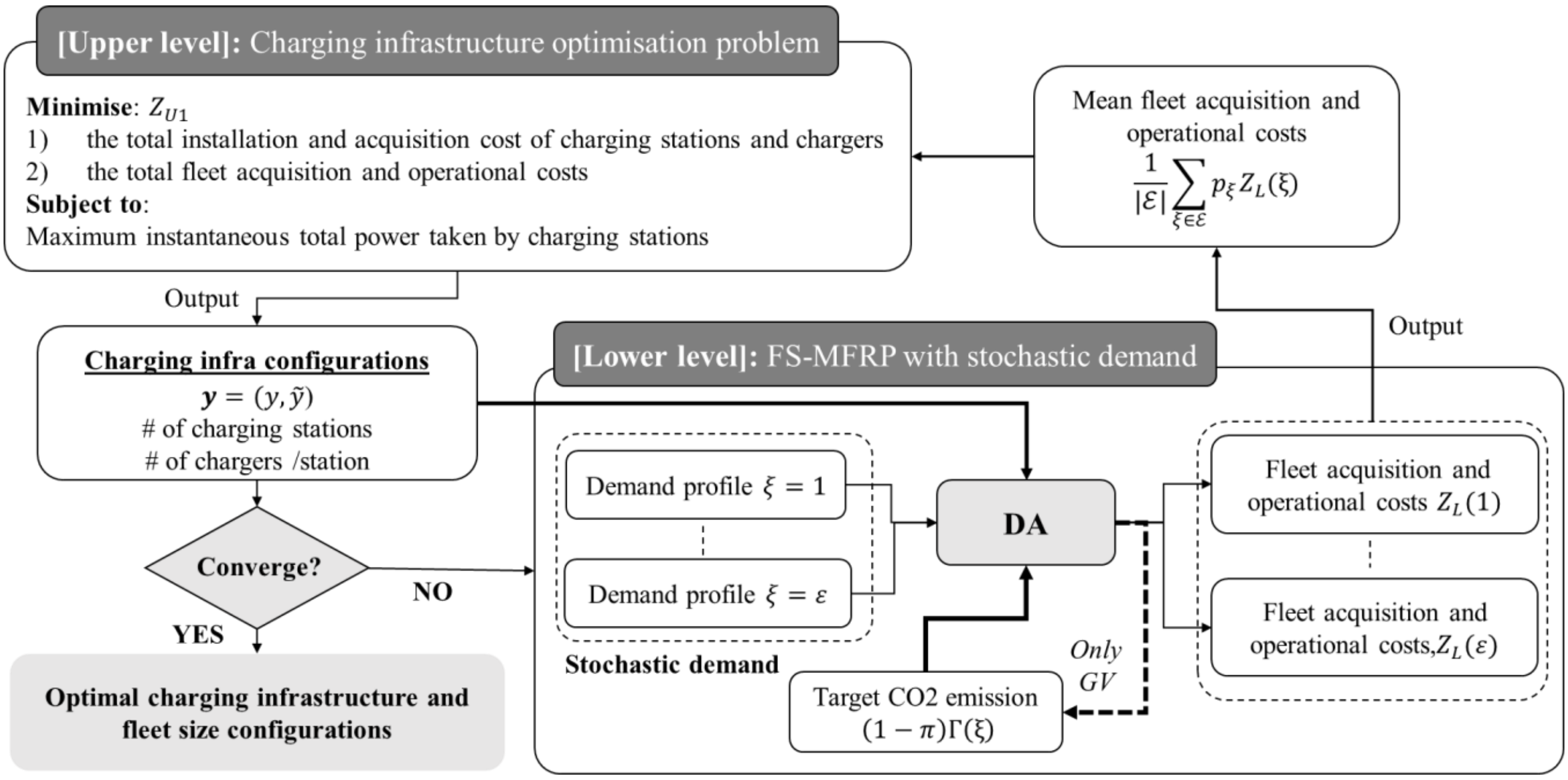

**Figure 1. The bi-level problem solution framework.**

**Algorithm 1. Bi-level charging infrastructure configuration and fleet size optimization algorithm.**

1. Input: Generate a set of demand scenarios $\mathcal{E}$ and initialize $W, \omega, y_{hw}^{max}$
2. Initialize a set of charging infrastructure configuration candidates $\tilde{S}$. Set $\tilde{s} = \tilde{s}_{best} = \emptyset$
3. Sort $\tilde{S}$ according to $n(\tilde{s}), \forall \tilde{s} \in \tilde{S}$
4. **for** $\tilde{s}$ in $\tilde{S}$
5.   **for** $\xi$ in $\mathcal{E}$
6.     Solve the FS-MFRP problem based on Algorithm 2 for each demand scenario
7.   **end for**
8.   Compute $\bar{Z}_L(\tilde{s})$
9.   **if** $\tilde{s}_{best} = \emptyset$
10.     $\tilde{s}_{best} \leftarrow \tilde{s}$
11.   **else**
12.     if $c(\tilde{s}) < c(\tilde{s}_{best})$ then $\tilde{s}_{best} \leftarrow \tilde{s}$
13.   **end**
14. **end for**
15. return $\tilde{s}_{best}$

Remark: $n(\tilde{s})$ is the number of chargers for the charging station configuration $\tilde{s}$. $\bar{Z}_L(\tilde{s})$ is the value of Eq. (6). $c(\tilde{s})$ is the value of Eq. (1).

To solve the FS-MFRP efficiently, we propose a deterministic annealing (**DA**) based metaheuristic (Braekers et al., 2014; Ma et al., 2024) with dedicated local search (**LS**) operators to minimize the total operational costs given a demand scenario. A penalty function approach is adopted to handle constraint violations for $CO_2$ emissions, vehicles' energy consumption, the number of used vehicles, and unserved customers. The unitary penalties associated with each component of the penalized cost function are tuned to have effective effects on the performance of the algorithm (described later). The sub-problem of charging scheduling under capacitated charging stations and partial recharge is solved by extending our previous work (Ma et al., 2024). Given a charging infrastructure configuration, a lower-level problem is solved twice for each demand scenario (i.e. once with only GVs to obtain the maximum $CO_2$ emission by dropping off Eq. (9), and then solve it again with Eq. (9) with a mixed fleet setting). The DA-based metaheuristic is described as follows (Algorithm 2). First, the problem instance is pre-processed to set up time windows associated with picked-up locations (for outbound requests) and drop-off locations (for inbound requests). Travel time (distance) of infeasible arcs (Sect. 3.2) is set with big positive values to forbid their appearance in the solution. The algorithm generates an initial solution $s_{init}$ as the best of 500 feasible solutions using a greedy insert approach (line 4 in Algorithm 2). This approach randomly selects an unserved request and inserts it into a random route at a feasible location with the least cost increase until all requests are served. The feasibility in terms of vehicle capacity, maximum riding time, and time window constraints is checked using the eight-step evaluation scheme (Cordeau and Laporte, 2003). The energy constraints (Eqs. (24)-(28)) are only checked for electric vehicles. When violated, a charging scheduling heuristic is applied to insert charging operations. Note that the charging synchronization constraints (Eqs. (29)-(33)) are checked once all unserved requests are inserted. The penalized cost function $\tilde{c}(s)$ of a solution $s$ is evaluated with four penalty terms as follows.

$$\tilde{c}(s) = c(s) + \rho_1 \tilde{n}_{cus}(s) + \rho_2 n_k(s) + \rho_3 max\big(0, \gamma(s) - \vec{\Gamma}\big) + \rho_4 \left( \sum_{k \in K^e} max\big(0, \vec{E}^k(s) - (E^k_{init} + \Delta E^k)\big) \right) \quad (40)$$

The first term $c(s)$ is the original objective function (8). The second term is the penalty for $\tilde{n}_{cus}(s)$ unserved users. The third term is the penalty for the number of used vehicles $n_k(s)$. The fourth term is the penalty for $CO_2$ constraint violation. Note that $\vec{\Gamma}$ denotes the targeted $CO_2$ emission amount of demand scenario $\xi$, i.e. $\vec{\Gamma} = (1-\pi)\Gamma(\xi)$. $\gamma(s)$ is the total $CO_2$ emission of solution $s$. The fifth term is the penalty when no feasible charging schedule can be found. This penalty is the summation of supplementary energy to be charged for vehicles to become energy-feasible to serve the customers. $\vec{E}^k(s)$ is the energy consumption of vehicle $k$ on solution $s$, and $E^k_{init}$ is the initial state of charge of vehicle $k$, and $\Delta E^k$ is the amount of energy can be feasibly charged for the current route without violating other constraints (time windows, vehicle capacity and maximum ride time of passengers). $\rho_1, \rho_2, \rho_3, \rho_4$ and other algorithmic parameters need to be tuned to achieve better performance (described later).

Following the generated initial solution $s_{init}$, a current solution $s$, a temporary solution $s'$, and a current best solution $s_{best}$ are initialized as $s_{init}$. Then $s$ and $s_{best}$ are updated iteratively based on a large neighbourhood search framework. An LS (local search) operator is randomly selected among a list of LS operators and applied for $s$ to obtain a temporary solution $s'$ (lines 11-12 in Algorithm 2). For every $\bar{n}_{remove}$ iterations, a *remove-route* operator is applied to reduce the number of used vehicles (lines 13-15 in Algorithm 2). If the penalized cost of $s'$ is lower than that of $s$ plus a threshold $T$ without charging operation conflicts for used EVs, the current solution $s$ is updated by $s'$ (lines 16-17 in Algorithm 2). When $s$ has no unserved customers, a $route_exchange\&improve$ operator is applied to exchange vehicle types for s given the $CO_2$ emission constraint. As the resulting solution might contain routes of EVs violating energy constraints, the $worst_relocate$ operator is applied on each energy-infeasible route. This operator removes the worst request (in terms of $\tilde{c}(s)$) of each energy-infeasible route and inserts them back to the least cost of other routes (lines 18-19 in Algorithm 2). If the resulting solution is promising, a repair procedure is evoked to repair its infeasibility related to charging conflicts and energy feasibility constraints.

Note that when charging facilities are very scarce compared to charging demand, finding a feasible charging plan for the fleet is impossible without updated charging station occupancy information. To address this situation, a $repair_sol$ procedure is evoked to find a feasible charging plan without charging conflicts between vehicles (lines 21-28). This procedure is composed of two steps. First, we check the energy feasibility of routes. If violated, randomly selected requests are ejected to the unserved pool until the energy feasibility is repaired. If the resulting solution has unserved customers or there are charging conflicts between vehicles, we repair it as follows. First, electric vehicles are sorted based on their charged energy in descending order. Initiate a charging state occupancy matrix (a binary matrix indicating whether a charger is occupied or not at time $t$). Then, iterate EVs with charging operations on the list and update the charging state occupancy matrix sequentially. If there are charging operation conflicts (charging time overlaps on chargers) with previous EVs, remove the charging operations of the current EVs and eject randomly selected requests to the unserved pool until this vehicle's route becomes energy-feasible. Continue this procedure until a non-conflict charging facility utilization is initialized. The resulting solution needs to re-insert the unserved requests on the list of EVs (no additional $CO_2$ emission). A greedy insertion policy is used to insert unserved requests one by one on k randomly selected (k=3) EVs. If there are unserved requests at the end, a new EV is used to insert them. Restart this procedure until all requests are served if there are still unserved requests.

After each request's insertion, vehicles' charging schedules (if any) are updated, given the updated charging state occupancy matrix information. This charging schedule variant adds an additional filter to the charging scheduling algorithm to ensure no charging conflicts between vehicles when adding charging operations in this solution repair process. Note that this $repair_sol$ procedure might be very computationally time-consuming if the charging supply is insufficient. To reduce unnecessary computational time, a Table list is maintained to avoid repairing the same infeasible solution previously visited (line 22). Note that an acceleration strategy can be used, allowing evoking this procedure with a predefined probability $\alpha$. Finally, if the cost or the number of used vehicles of the resulting feasible solution $s$ is smaller than those of $s_{best}$ and the solution is feasible, $s_{best}$ is updated by $s$ (lines 30-31). If the non-improvement count $i_{imp}$ is greater than 0 (i.e. $s_{best}$ remains unimproved during the current iteration), then the threshold $T$ is reduced by $T_{max}/T_{red}$. If $T$ becomes negative, reset $T$ randomly between 0 and $T_{max}$ (line 38 in Algorithm 2). If $i_{imp} > n_{imp} * n_k(\boldsymbol{s_{best}})$**,** where $n_k(s_{best})$ is the number of vehicles of $s_{best}$, reset the temporary solution $s'$ as $s_{best}$ (restart) and $i_{imp} = 0$ (lines 39 and 40 in Algorithm 2). Note that the feasibility checks for charging conflicts, $CO_2$ violation and energy feasibility are implemented in an efficient way (time complexity is $O(n)$ or so, where $n$ is the number of EVs) using information stored during the current solution search process.

**Algorithm 2. DA-based metaheuristic algorithm for solving the lower-level FS-MFRP**.

1. Input: customer requests, charging stations, and the parameters of the system
2. Pre-processing (time window tightening and infeasible arcs eliminations)
3. Set up e-DARP instance
4. $s_{init} \leftarrow$***generateInitSol****(e-DARP)*
5. Set the current solution $s = s' = s_{best} = s_{init}$, where $s'$ is a temporary solution. $i_{imp} = 0$, $T = T_{max}$
6. **For** $iter = 1: iter_{max}$
7. $\quad i_{imp} \leftarrow i_{imp} + 1$
8. $\quad$ **if** $count_{stagnant} = \bar{n}_{stagnant}$
9. $\quad\quad$ return $s_{best}$
10. $\quad$ **end if**
11. $\quad ls \leftarrow rand(LS)$
12. $\quad s' \leftarrow ls(s')$
13. $\quad$ If $iter$ % $\bar{n}_{remove} = 0$
14. $\quad\quad s' \leftarrow remove_route(s')$
15. $\quad$ end
16. $\quad$ **if** $(\tilde{c}(s') < \tilde{c}(s) + T$
17. $\quad\quad s \leftarrow s'$
18. $\quad\quad$ **if** $n_{\text{unserved}}(s) = 0$
19. $\quad\quad\quad s \leftarrow route_exchange\&improve(s)$
20. $\quad\quad\quad$ **if** $(c(s) < c(s_{best}) || \, n_k(s) < n_k(s_{best}))$ && $CO_2$ emission is feasible
21. $\quad\quad\quad\quad$ If not $no_charging_conflict(s)$ *or not energy_feasible(s)*
22. $\quad\quad\quad\quad\quad$ If not $tabu_list(s)$ && $random(0,1) < \alpha$
23. $\quad\quad\quad\quad\quad\quad repair_energy_feasibility(s)$
24. $\quad\quad\quad\quad\quad\quad$ If not $n_{\text{unserved}}(s) = 0$
25. $\quad\quad\quad\quad\quad\quad\quad repair_sol$(s)
26. $\quad\quad\quad\quad\quad\quad$ end
27. $\quad\quad\quad\quad\quad$ end
28. $\quad\quad\quad\quad$ end
29. $\quad\quad\quad$ end
30. $\quad\quad\quad$ **if** $(c(s) < c(s_{best})$ && $no_charging_conflict(s)$ && $energy_feasible(s)$

31. $s_{best} \leftarrow s; i_{imp} = 0$
32. **end if**
33. **end if**
34. **end if**
35. **if** $i_{imp} > 0$
36. $T \leftarrow T - T_{max}/T_{red}$
37. **if** $T < 0$
38. $T \leftarrow random(0,1) * T_{max}$
39. **if** $i_{imp} > n_{imp} * n_k(\boldsymbol{s}_{best})$
40. $s' \leftarrow s_{best}; i_{imp} = 0$
41. **end if**
42. **end if**
43. **end if**
44. **end for**
45. return $s_{best}$

We propose a list of customised LS operators to find feasible and improved solutions in the neighbourhood of the current solution $s$ as follows. Four new solution improvement strategies are added to address the $CO_2$-emission-constrained mixed fleet size and routing problem: i) a remove-route operator is proposed to minimise the number of used vehicles; ii) a tailored penalized cost function (Eq. (40)) is proposed for penalizing violated constraints to enlarge the searched neighbourhood from a temporal solution; iii) a route-exchange operator is proposed to swap vehicle types given the target $CO_2$ constraints; iv) We propose an infeasibility-repair strategy to improve the chance of finding feasible charging schedules in a congested charging station situation when the number of fast chargers is very limited.

Note that the energy constraints are checked at the end of each LS operator to determine whether vehicles need to go recharge or not. The charging synchronization constraint is only checked for promising solutions (lines 21 in Algorithm 2). The implementation technique for the solution representation is based on the double-linked array data structure (Kytöjoki et al., 2007), while the charging operations are managed using another array to store the sequence of charging operations, each of which includes the predecessor, the charger's ID, and the start and end times of the charging operations. The proposed LS operators are described as follows.

- **Relocate ensemble**: Randomly select one of the following two relocate operators: greedy relocate and worst relocate. The first randomly removes a request from one route and inserts it back to the least-cost feasible position of the other routes. The second removes the most expensive request on a randomly selected route and inserts it back into the feasible and least-cost position of the other routes.
- **Destroy-repair**: This operator removes (destroys) a part of the current solution (served requests) using a set of remove operators and then repairs them with two insertion operators. The degree of destruction (pre-defined) is controlled by a parameter $\delta \in [0.2, 0.5]$, representing the percentage of the current solution to be destroyed. First, a randomly selected removal operator removes served requests to a pool of unserved users. The removal operators include random-removal, worst-removal, distance-related removal, time-window-related removal, and Shaw-removal. Following the removal procedure, a randomly selected repair operator is applied to insert back unserved requests into current routes. At the end of this procedure, if there are still unserved requests, create a new route to insert these unserved requests. Remaining unserved requests are kept in the unserved pool. Note that when creating a new route, we need to determine its vehicle type (electric or gasoline). This decision is determined automatically according to whether the current solution (before applying the destroy-repair operator) violates $CO_2$ constraints or not (i.e., if violated, create a new EV; otherwise, create a gasoline vehicle). Two repair operators are considered: greedy insert and regret insert. The first inserts back unserved requests to a feasible position with the least penalized costs of the other routes. If the inserted route (vehicle) is an EV that needs to be recharged, apply the charging schedule algorithm to schedule charging operations. Unlike our previous study, we allow vehicles' energy constraints to be violated with a penalty cost when no feasible charging operations can be scheduled. For the regret insert operator (Ma et al., 2024), we insert back unserved users based on k-regret strategies (k=2 or 3). The penalized cost function is defined as Eq. (40).

- **Two-opt**: Determine the length of the to-be-operated segment as a random integer between 2 and 4. Then reverse the visit order of the segment along the route until a feasible and improved solution is found.
- **Four-opt**: Remove four consecutive arcs (three nodes) and find the best visit order for that segment among all possible permutations. Continue this operation along the route and retain the feasible solution with the least cost.
- **Two-opt***: Randomly select two routes. Identify a list of candidate arcs where the vehicle load is zero. For each route, randomly remove one candidate arc and combine the first part of the first route with the second part of the second route along the route until a feasible and improved solution is found.
- **Swap-requests**: Swap two requests on two different routes. The position of pick-up (drop-off) of the first request needs to be reinserted at the position of pick-up (drop-off) of the second request. If feasible, the removed request of the second route is reinserted to a feasible position with the least cost on the first route. If failed, the request is reinserted to another route until a feasible and improved solution is found.
- **Swap-segments**: Randomly select two routes and identify the segments of the routes where the vehicle load is zero. Exchange the segments along the route and retain the feasible solution with the least cost.
- **Remove-route**: Let $R$ be the set of all requests. We denote $R^{+}(s)$ and $R^{-}(s)$ as the subsets of served and unserved requests of the current solution s, respectively. Randomly select a route and put their requests on $R^{-}(s)$. Then remove $m$ requests ($m$ is a random integer between 1 and $min(n, \delta|R^{+}(s)|)$, where $n$ is the pre-defined maximum number of requests can be removed from $R^{+}(s)$, $\delta$ is the degree of destruction between 0.2 and 0.5, and $|R^{+}(s)|$ is the number of served requests on $s$. Randomly select one removal operator and apply it to remove $m$ served requests and add them to $R^{-}(s)$. Then use the greedy insert or regret insert operator to insert unserved requests in $R^{-}(s)$ to used vehicles or one newly created vehicle.

The different roles and control parameters of the operators of the DA metaheuristic are described in Table 1.

Table 1. The different roles and control parameters of the operators of the DA metaheuristic.

| Operator | Type | Description | Involved route | Control parameter |
|---|---|---|---|---|
| Relocate ensemble | LS | RC reduction | inter-route | Randomizing the use of the greedy and worst relocate operator |
| Destroy-repair | LS | Shaking | inter-route | Using the degree of destruction to determine the number of removed customers |
| Random-removal | Auxiliary | Shaking | inter-route | |
| Worst-removal | Auxiliary | Shaking | inter-route | Degree of randomness parameter |
| Distance-removal | Auxiliary | Shaking | inter-route | |
| TW-removal | Auxiliary | Shaking | inter-route | |
| Shaw-removal | Auxiliary | Shaking | inter-route | Relatedness parameters |
| Two-opt | LS | RC reduction | *Intra-route* | |
| Four-opt | LS | RC reduction | *Intra-route* | |
| Two-opt* | LS | RC reduction | inter-route | |
| Swap-requests | LS | RC reduction | inter-route | |
| Swap-segments | LS | RC reduction | inter-route | |
| Remove-route | Global control | Fleet size reduction | inter-route | |
| Route_exchange & improve | Global control | Fleet composition | inter-route | |
| Greedy relocate | Auxiliary | RC reduction | inter-route | |
| Worst relocate | Auxiliary | RC reduction | inter-route | |
| Greedy insert | Auxiliary | RC reduction | inter-route | |
| Regret insert | Auxiliary | RC reduction | inter-route | k-regret parameter |

Note: LS: Local search; RC: Routing costs; TW: Time-window

## 4. Numerical study

In this section, two numerical experiments are described; we first test the methodology on small test instances to tune algorithmic parameters and validate the DA algorithm to solve the lower-level FS-MFRP problem, and then we demonstrate the applicability of the proposed models and algorithms through a realistic case study focusing on the feeder service area around the Bettembourg train station in Luxembourg. Bettembourg is surrounded by a rural hinterland with dispersed villages, making feeder services highly relevant.

### 4.1. Model parameter and test instance generation

We consider a synthetic stylised test network consisting of a 256km square service area (with the side length of 16km ) with one train station located at the centre and one depot equipped with overnight chargers. The train is assumed to be operated from 6 a.m. to 11 p.m. with varying frequencies; 15 minutes between 7-9 a.m. and 5-7 p.m., 30 minutes between 6-7 a.m., 9-11 a.m., 4-5 p.m. and 7-8 p.m., and 1 hour during the rest of service time. DRFS provides first- and last-mile services to/from the train station, where users can be picked up/dropped off at potential stops or meeting points. The potential stops are uniformly distributed in the service area at a 1 km distance apart to randomly generate demand scenarios. The demand is defined in terms of service requests, each request containing one to four passengers. The maximum user acceptable waiting time at the station is set as 10 minutes, and the maximum in-vehicle time is also set as 1.5 times the direct travel time between customers' pickup and drop-off locations to account for a maximum acceptable detour time for the passengers. Service time at each pickup location is assumed to be 30 seconds.

For the first experiment, we assume that there are 2 DC fast chargers (50 kWh) installed at one charging station located at the transit station (see Figure 2). The considered vehicle types are the type-I EV and GV with the same passenger capacity (24 seats). The problem is to determine the mixed fleet size to minimize total operational costs under the target $CO_2$ emission constraint. The experiment setting for the Bettembourg case study will be described later. The characteristics of vehicles and charging stations are summarized in Table 2.

The proposed solution algorithms are programmed in Julia and run on a laptop with Intel(R)Core (TM) i7-11850H processor and 32 GB memory using a single thread. We use Gurobi v10.0.3 to solve the mixed integer program of FS-MFRP to obtain the exact solutions.

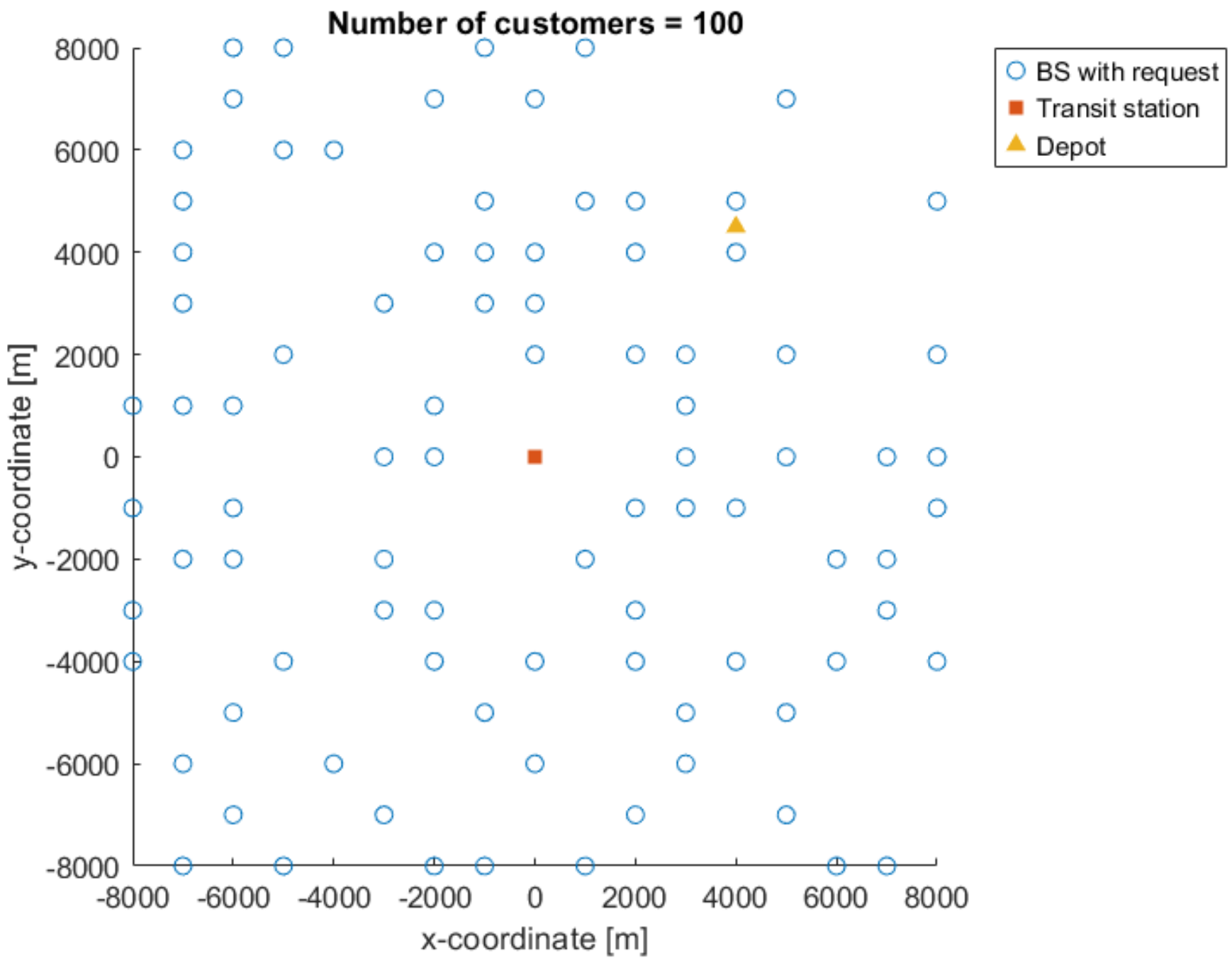

**Figure 2. Example of the test instance with 100 requests with one single depot, transit station and charging station.**

**Table 2. The characteristics and costs of different types of vehicles and charging infrastructure.**

| **Gasoline/Electric vehicles** | **EV[1]** | | **GV[1]** |
|---|---|---|---|
| | **Type 1** | **Type 2** | |
| Vehicle Capacity (number of seats) | 24 | 10 | 24 |
| Fuel Capacity (kWh) or (l) | 117 | 85 | 10000[2] |
| Energy consumption rate (kWh/km) or (l/km) | 0.938 | 0.469 | 0.067 |
| Energy cost (€/kWh) or (€/l) | 0.23 | 0.23 | 1.83 |
| Daily equivalent of purchasing and maintenance cost (€/day) | 23.78 | 11.89 | 16.17 |
| Average speed (km/h)[3] | 50 | 50 | 50 |
| CO2 emission (kg/km) | 0 | 0 | 0.176 |
| ***Charging infrastructure*** | **Rapid** | | **Super-fast** |
| Charging station opening cost (€/day) | 4.11 | | 68.08 |
| Cost of rapid charger (€/day) | 9.59 | | 19.18 |
| Charging power (kW/hour) | 125 | | 220 |

Remark:1. Adapted from Meishner and Sauer (2020). 2. We assume that GVs can get refuelled at nearby gasoline stations with negligible access time, equivalent to an unlimited fuel capacity.

4.2. Algorithmic parameter tuning

The DA algorithm for the FS-MFRP needs to tune its algorithmic parameters to be effective to solve the lower-level mixed fleet size and routing problem (there are 10 parameters to be tuned). We generate 5 test instances with 20, 40, 60, 80, and 100 requests, randomly sampled from the potential demands associated with vehicle stops as aforementioned. The initial battery level of vehicles is set as 50%-full. We consider a 50% $CO_2$ emission reduction target for the parameter tuning. To tune these parameters efficiently, we adopt a sequential approach procedure by tuning a good reference parameter setting. Note that more advanced parameter-tuning approaches can also be adopted (Bischl, et al., 2023). First, we conduct preliminary experiments and identify a reference parameter set,
i.e. $(\rho_1, \rho_2, \rho_3, \rho_4, \bar{n}_{remove}, \bar{n}_{stagnant}, iter_{max}, n_{imp}, T_{red}, t_{max})$= (40, 100, 30, 1, 250, 150, 10^5, 300, 300, 0.6). Then we set five test values for each parameter. For each time, we evaluate the performance of the algorithm by changing the value of one parameter with respect to the reference parameter set. The best values for previously tested parameters are retained and used for the following experiments to tune the remaining parameters. In this way, we can obtain a good parameter setting with reasonable computational time if compared with exhaustive testing for all possible permutations of test parameter values. Each test instance is solved 3 times to get the average objective function value. The tested parameter values are set as follows.

- $\rho_1$: Penalty for each unserved user (Eq. (40)), $\rho_1 \in \{20,40,60,80,100\}$
- $\rho_2$: Penalty for each used vehicle (Eq. (40)), $\rho_2 \in \{100,150,200,250,300\}$
- $\rho_3$: Penalty for excess $CO_2$ emission compared to the targeted threshold (Eq. (40)), $\rho_3 \in \{10,20,30,40,50\}$
- $\rho_4$: Weighted coefficient for each kWh charged to satisfy the energy constraint (Eq. (40)), $\rho_4 \in \{0.1, 0.3, 0.5, 0.7, 0.9\}$
- $\bar{n}_{remove}$: Run $remove_route$ LS operator every $\bar{n}_{remove}$ iterations (line 13 in Algorithm 2), $\bar{n}_{remove} \in \{100,150,200,250,300\}$
- $\bar{n}_{stagnant}$: Maximum number of stagnation (line 8 in Algorithm 2) after which the algorithm is stopped, $n_{stagnant} \in \{50,100,150,200,250\}$
- $iter_{max}$: Maximum number of iterations (line 6 in Algorithm 2), $iter_{max} \in \{50\text{k}, 100\text{k}, 150\text{k}, 200\text{k}, 300\text{k}\}$
- $n_{imp}$: Restart parameter (line 30 in Algorithm 2), $n_{imp}$ ∈{100, 200, 300, 400, 500}

- $T_{red}$ : Threshold reduction factor for reducing the threshold value ($T \coloneqq T - \frac{T_{max}}{T_{red}}$) (line 27 in Algorithm 2), $T_{red} \in \{100, 200, 300, 400, 500\}$
- $T_{max}$: A user-defined coefficient to determine the maximum threshold value $T_{max}$ to accept worsen temporary solutions (line 27 in Algorithm 2), i.e. $t_{max} \in \{0.3, 0.6, 0.9, 1.2, 1.5\}$.

Table 3 reports the average objective function values obtained by the DA algorithm over the 5 test instances for each tested parameter value. The bold texts are the best parameter values for their respective parameters. We can observe the impact of different parameter values on the quality of obtained solutions. It is notable that many of these results show that the behaviour is not simply quadratic. The objective function value does not always get monotonically higher (worse) when going further away from the chosen value. The retained best parameter setting, is then used for the numerical studies in the following sections i.e., $(\rho_1, \rho_2, \rho_3, \rho_4, \bar{n}_{remove}, \bar{n}_{stagnant}, iter_{max}, n_{imp}, T_{red}, T_{max}) = (100,200,20, 0.3,250,150,10^5, 300,400,0.9)$.

**Table 3. Sequential parameter tuning for the DA algorithm.**

| Parameter | Tested value | | | | |
|---|---|---|---|---|---|
| $\rho_1$ | 20 | 40 | 60 | 80 | **100** |
| Avg. obj. value | 170.71 | 170.31 | 169.31 | 167.49 | 166.47 |
| $\rho_2$ | 100 | 150 | **200** | 250 | 300 |
| Avg. obj. value | 169.91 | 170.93 | 167.09 | 187.10 | 177.96 |
| $\rho_3$ | 10 | **20** | 30 | 40 | 50 |
| Avg. obj. value | 177.47 | 167.83 | 181.72 | 170.10 | 171.73 |
| $\rho_4$ | 0.1 | **0.3** | 0.5 | 0.7 | 0.9 |
| Avg. obj. value | 170.60 | 166.87 | 188.11 | 172.44 | 167.36 |
| $\bar{n}_{remove}$ | 100 | 150 | 200 | **250** | 300 |
| Avg. obj. value | 173.64 | 179.81 | 167.46 | 167.35 | 178.62 |
| $\bar{n}_{stagnant}$ | 50 | 100 | **150** | 200 | 250 |
| Avg. obj. value | 182.27 | 197.94 | 166.31 | 168.65 | 170.01 |
| $iter_{max}$ | 50K | **100k** | 150k | 200k | 300k |
| Avg. obj. value | 175.73 | 166.12 | 167.55 | 170.14 | 188.19 |
| $n_{imp}$ | 100 | 200 | **300** | 400 | 500 |
| Avg. obj. value | 174.86 | 167.24 | 165.85 | 169.01 | 181.33 |
| $T_{red}$ | 100 | 200 | 300 | **400** | 500 |
| Avg. obj. value | 170.71 | 178.59 | 191.55 | 169.70 | 175.66 |
| $T_{max}$ | 0.3 | 0.6 | **0.9** | 1.2 | 1.5 |
| Avg. obj. value | 179.39 | 196.86 | 167.56 | 173.30 | 168.98 |

Remark: Average objective function value is based on the 3-run average over 5 test instances (c20, c40, c60, c80, c100) with initial battery level of 50%-full and 50% $CO_2$ emission reduction target.

4.3. Computational results to validate the DA algorithm

We test the performance of the DA algorithm for the lower-level problem with 0%, 50% and 90% $CO_2$ reduction targets and compare it with the solutions obtained by a state-of-the-art MIP solver with a 4-hour computational time limit. Given that the commercial solver can obtain exact solutions up to a small problem size, we limit ourselves to the problem size with up to 50 requests considering the computational time limit. In total, 15 instances (each problem size, i.e. 10, 20,..,50 requests, has three test instances) are tested with the number of requests ranging from 10 to 50, each request containing 1 to 4 passengers. The initial battery levels of vehicles are kept at 50% full, similar to the experimental settings used for parameter tuning. The parameter setting is shown in Table 2. The solutions obtained by the commercial solver are shown on the block of MILP columns. We can observe that for the gasoline fleet scenario, the commercial solver can obtain exact or near-optimal solutions for all the test instances. When reducing the $CO_2$ emission target, the problem becomes more difficult to solve exactly.

The performance of the DA algorithm is shown in Table 4, on the right block of DA using the reference parameter setting obtained in Section 4.2. For the scenario of all GVs, the best and average gaps to the best-known solutions (BKS) are almost 0 for all test instances. When reducing 50% $CO_2$ emission, the best and average gap to BKS is no more than 2.23% and 4.26% for all tested instances, respectively. For the 90% $CO_2$ reduction scenario, the best and average gap to BKS is no more than 4.10% and 6.16% for all tested instances, respectively. The computational time for the DA algorithm is no more than 7 seconds for the worst cases. The $CO_2$ emissions of the obtained solutions satisfy the target levels. The result demonstrates that the proposed DA algorithm can solve the mixed fleet size and routing problem with different levels of $CO_2$ emission constraints efficiently with good solution quality.

To evaluate the performance of the DA algorithm for larger instances, we conducted experiments with the number of requests set to 100, 200, 300, 400, and 500 while maintaining consistent parameter settings. Table 5 summarizes the results for three scenarios: (1) GV-only fleet, (2) mixed fleet with fully charged EVs, and (3) mixed fleet with initial battery levels set to 50%. For scenarios (2) and (3), two $CO_2$ reduction targets of 50% and 90% were considered.

The results indicate that when EVs are fully charged, no charging behaviour is observed, even in the case of 500 requests. This is because the fleet size required to fulfil all requests is sufficiently large, resulting in shorter vehicle kilometre travel per vehicle, thereby eliminating the need for charging when the battery is initially fully charged. This observation is further supported by the fact that the fleet sizes in the GV-only and mixed fleet cases are nearly identical, suggesting that the mixed fleet with fully charged batteries requires only a minimal number of vehicles to meet the service demands. This is the case for the mixed fleet with an initial battery level of 50%. In this scenario, charging behaviour is observed, as the initial battery level is insufficient to complete the service without recharging. Incorporating charging behaviour is undoubtedly more cost-effective than purchasing additional EVs, demonstrating that the DA algorithm performs as expected.

Given the significant computational time required to solve these large instances using MILP, we did not include it as a reference case. Consequently, it was not possible to quantify the gaps between the DA solutions and the exact solutions. Nevertheless, the results show that the DA algorithm successfully produces objective values while meeting $CO_2$ emission constraints, utilizing a mixed fleet, and incorporating charging activities (see Table 5).

**Table 4. The performance of the DA-based metaheuristics.**

| | | MILP | | | | | DA | | | | |
|---|---|---|---|---|---|---|---|---|---|---|---|
| CO2 emission target* | # requests | Best known obj. (BKS) ** | Gap to the lower bound | $CO_2$ (Best sol.) (kg) | Total charging time (min.) | cpu (sec.) | Best Gap to BKS (%) | Avg. Gap to BKS (%) | $CO_2$ (Best sol.) (kg) | Total char. time (min.) | cpu (sec.) |
| 10% | c10 | 60.71 | 0.00% | 0.0 | 0.0 | 8 | 0.00 | 0.00 | 1.33 | 0.00 | 1 |
| | c20 | 92.85 | 4.08% | 1.0 | 6.5 | 7263 | 3.23 | 3.26 | 2.00 | 0.00 | 1 |
| | c30 | 130.62 | 15.20% | 1.8 | 0.0 | 14400 | 4.10 | 6.16 | 3.00 | 0.33 | 5 |
| | c40 | - | - | - | - | 14400 | na | na | 3.33 | na | 5 |
| | c50 | - | - | - | - | 14400 | na | na | 4.00 | na | 7 |
| 50% | c10 | 54.36 | 0.00% | 3.7 | 0.0 | 13 | 0.00 | 0.00 | 1.33 | 0.00 | 1 |
| | c20 | 67.13 | 0.00% | 19.4 | 0.0 | 5909 | 0.29 | 1.08 | 2.00 | 0.00 | 1 |
| | c30 | 90.53 | 0.88% | 29.1 | 0.0 | 11358 | 2.13 | 4.26 | 2.67 | 0.00 | 3 |
| | c40 | 116.51 | 23.62% | 36.8 | 19.2 | 14400 | 2.23 | 3.94 | 3.00 | -0.33 | 5 |
| | c50 | - | - | - | - | 14400 | na | na | 3.00 | na | 6 |
| GV | c10 | 25.20 | 0.00% | 23.6 | 0.0 | 5 | 0.00 | 0.00 | 1.33 | 0.00 | 1 |
| | c20 | 38.33 | 0.00% | 39.0 | 0.0 | 78 | 0.00 | 0.01 | 2.00 | 0.00 | 0 |
| | c30 | 52.10 | 0.10% | 58.3 | 0.0 | 4901 | 0.00 | 0.04 | 2.67 | 0.00 | 1 |
| | c40 | 59.90 | 0.15% | 74.0 | 0.0 | 9882 | 0.05 | 0.17 | 3.00 | 0.00 | 1 |
| | c50 | 57.14 | 0.10% | 91.1 | 0.0 | 13513 | 0.14 | 0.22 | 2.67 | 0.00 | 1 |

Remark: The results are based on 5 runs, averaged over three datasets, each containing 5 randomly generated test instances with the number of requests ranging from 10 to 50 requests.

**Table 5. The performance of DA algorithm with large number of requests.**

| Init. SOC | CO2 target | # of requests | Best Obj. | Ave. Obj. | Fleet size GV | Fleet size EV | Fleet size Total | $CO_2$ best | Max $CO_2$ constr. | Char. time | cpu (sec.) |
|---|---|---|---|---|---|---|---|---|---|---|---|
| Full | GV only | 100 | 148.62 | 149.07 | 8 | 0 | 8 | 125.17 | - | 0 | 11 |
| | | 200 | 259.15 | 266.79 | 14 | 0 | 14 | 212.93 | - | 0 | 22 |
| | | 300 | 334.44 | 338.56 | 18 | 0 | 18 | 281.89 | - | 0 | 34 |
| | | 400 | 410.66 | 427.70 | 22 | 0 | 22 | 356.91 | - | 0 | 60 |
| | | 500 | 534.91 | 545.11 | 29 | 0 | 29 | 428.73 | - | 0 | 89 |
| | 0.5 | 100 | 238.47 | 251.66 | 5 | 3 | 8 | 62.18 | 62.59 | 0 | 14 |
| | | 200 | 434.95 | 447.49 | 8 | 6 | 14 | 105.95 | 106.47 | 0 | 26 |
| | | 300 | 555.42 | 566.10 | 11 | 7 | 18 | 140.94 | 140.94 | 0 | 40 |
| | | 400 | 710.09 | 722.34 | 14 | 10 | 24 | 176.55 | 178.45 | 0 | 52 |
| | | 500 | 864.83 | 876.75 | 17 | 12 | 29 | 214.04 | 214.36 | 0 | 77 |
| | 0.9 | 100 | 320.36 | 325.93 | 1 | 7 | 8 | 11.57 | 12.52 | 0 | 13 |
| | | 200 | 558.03 | 572.70 | 2 | 12 | 14 | 21.21 | 21.29 | 0 | 31 |
| | | 300 | 698.45 | 723.86 | 3 | 14 | 17 | 28.17 | 28.19 | 0 | 42 |
| | | 400 | 912.82 | 928.98 | 3 | 20 | 23 | 34.61 | 35.69 | 0 | 75 |
| | | 500 | 1162.08 | 1167.79 | 5 | 25 | 30 | 42.81 | 42.87 | 0 | 93 |
| 50% | 0.5 | 100 | 250.70 | 256.24 | 4 | 4 | 8 | 61.75 | 62.59 | 13.38 | 26 |
| | | 200 | 448.96 | 464.70 | 7 | 7 | 14 | 106.37 | 106.47 | 31.09 | 44 |
| | | 300 | 582.90 | 593.99 | 10 | 9 | 19 | 137.42 | 140.94 | 14.66 | 71 |
| | | 400 | 749.98 | 760.60 | 13 | 12 | 25 | 173.43 | 178.45 | 18.60 | 128 |

| | | | | | | | | | | |
|---|---|---|---|---|---|---|---|---|---|---|
| | 500 | 886.21 | 909.24 | 15 | 14 | 29 | 212.10 | 214.36 | 35.87 | 184 |
| 0.9 | 100 | 327.07 | 330.06 | 1 | 7 | 8 | 12.49 | 12.52 | 29.17 | 27 |
| | 200 | 584.08 | 592.16 | 3 | 12 | 15 | 20.18 | 21.29 | 30.31 | 45 |
| | 300 | 770.17 | 784.41 | 2 | 17 | 19 | 25.23 | 28.19 | 8.26 | 83 |
| | 400 | 975.07 | 989.15 | 3 | 21 | 24 | 35.04 | 35.69 | 47.86 | 130 |
| | 500 | 1176.72 | 1206.00 | 4 | 26 | 30 | 42.08 | 42.87 | 24.93 | 145 |

Remark: The results are based on 5 randomly generated test instances with the number of requests from 100 to 500 requests.

## 4.4. Case study

### 4.4.1. Case study settings and service characteristics

We consider a case study that mimics a real-world scenario for commuters using Bettembourg station as a transfer hub in their journey. Bettembourg station is located at the centre of the municipality of Bettembourg in Luxembourg and acts as a transit hub for commuters in Luxembourg and cross-border commuters from France. There are 3 train lines and 8 bus lines connected to this station. In 2022, on average, 3,756 people used the service, and 11,325 people got on/off the trains per day at this station (Digital Mobility Observatory, 2023). Currently, there are two companies providing on-demand feeder service around Bettembourg; however, but their operations are limited to a very small scale[4,5]. For this study, we assume that an existing operator of on-demand feeder services with GVs is planning to extend their service to the Bettembourg station. To meet the environmental goals, they would like to electrify (part of) their fleet upon the start of the service. Specification of the service, service area, and assumed demand characteristics are described below.

The service area is assumed to be a circle with a diameter of 17.8 km around Bettembourg station. Figure 3 illustrates the main train stations and bus stops connected to Bettembourg station within the service area. Luxembourg station (at the top) and Esch-sur-Alzette station (at the left bottom side in Figure 3) are the two biggest cities in Luxembourg and, hence, the two most popular stations, where 57,962 and 14,372 people get on/off the bus and train per day, respectively. The service area includes Esch-sur-Alzette station but not Luxembourg station. There are direct buses and trains connecting Esch-sur-Alzette station and Luxembourg station. However, depending on the departure time, the commuter sometimes needs to change trains/buses at Bettembourg station to travel between those two big stations.

The service is assumed to provide the first- and the last-mile service to/from the Bettembourg train station. Users whose origins/destinations are within 1 km from the Bettembourg station cannot use this service as they can walk to the station. The potential customers are 1) public transit users who make a transfer at Bettembourg station to get on/off the train and have their origin or destination in the service area (see Figure 4), and 2) private car users who travel from/to somewhere in a service area to/from the place reachable by a train line passing by Bettembourg station (see Figure 5). As we do not investigate the modal shift caused by introducing DRFS, we do not distinguish these two types of customers.

[4] Sales Lentz

[5] Services by Beetebuerg (e.g., Proxibus - Beetebuerg and eisen e-bus - Beetebuerg)

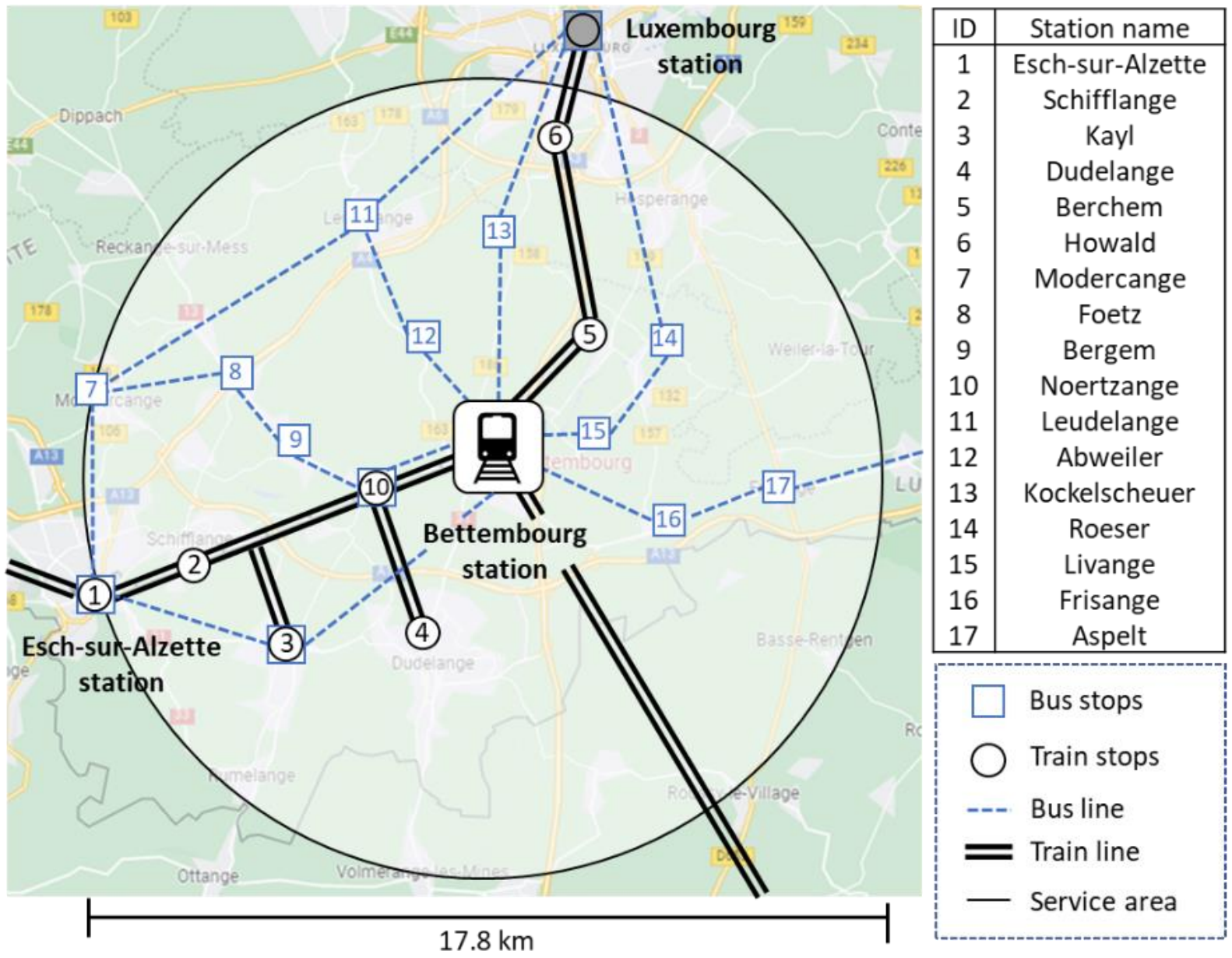

| ID | Station name |
|---|---|
| 1 | Esch-sur-Alzette |
| 2 | Schifflange |
| 3 | Kayl |
| 4 | Dudelange |
| 5 | Berchem |
| 6 | Howald |
| 7 | Modercange |
| 8 | Foetz |
| 9 | Bergem |
| 10 | Noertzange |
| 11 | Leudelange |
| 12 | Abweiler |
| 13 | Kockelscheuer |
| 14 | Roeser |
| 15 | Livange |
| 16 | Frisange |
| 17 | Aspelt |

**Figure 3. Service area assumed for this case study. The left top table displays the list of station names corresponding to the number of bus/train stops on the map.**

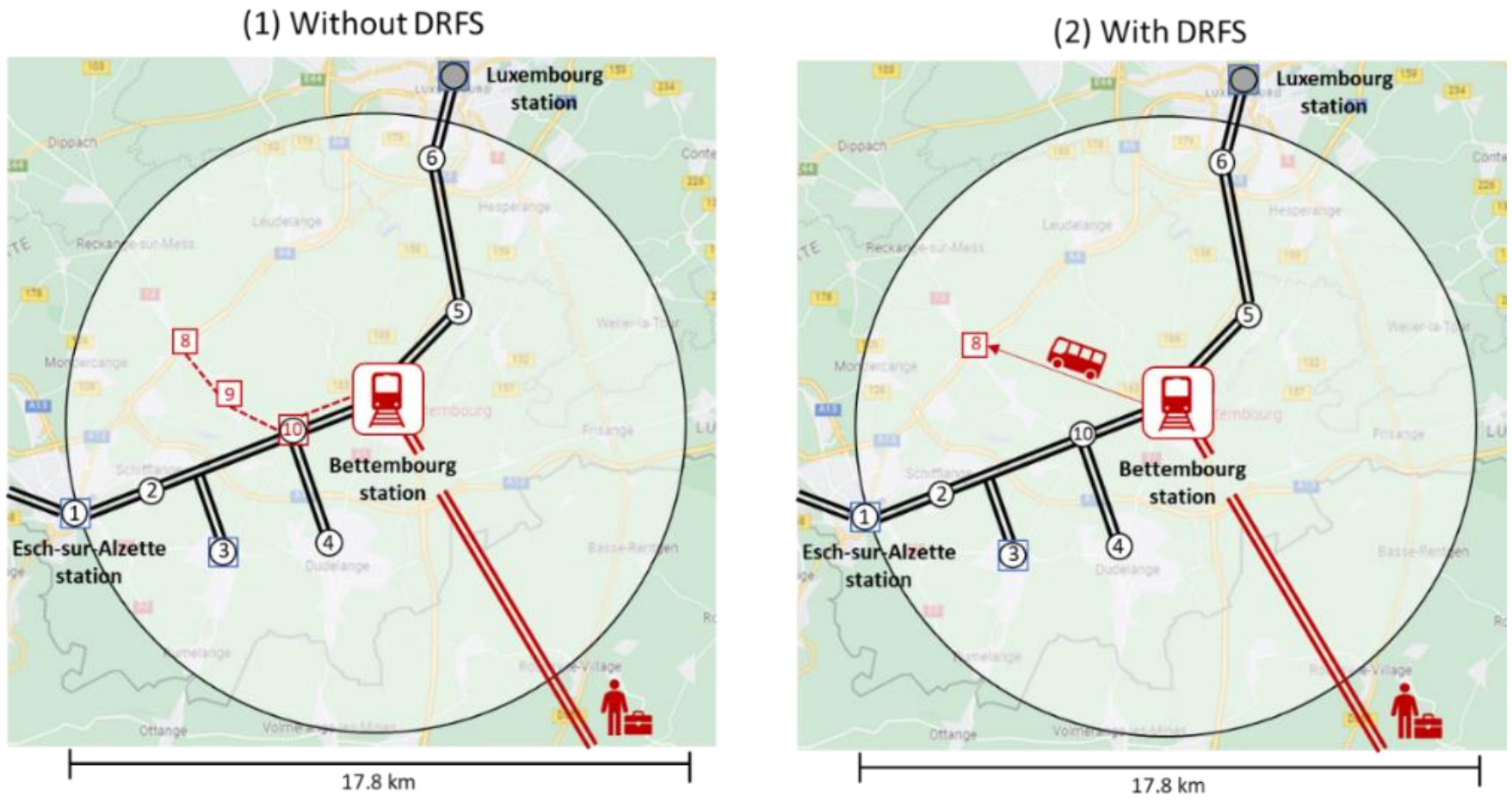

**Figure 4. Example of last mile services (the right figure) replacing the bus travel (the left figure) from the Bettembourg station to Foetz bus station (8) where a traveller is coming from France to the Bettembourg station by train.**

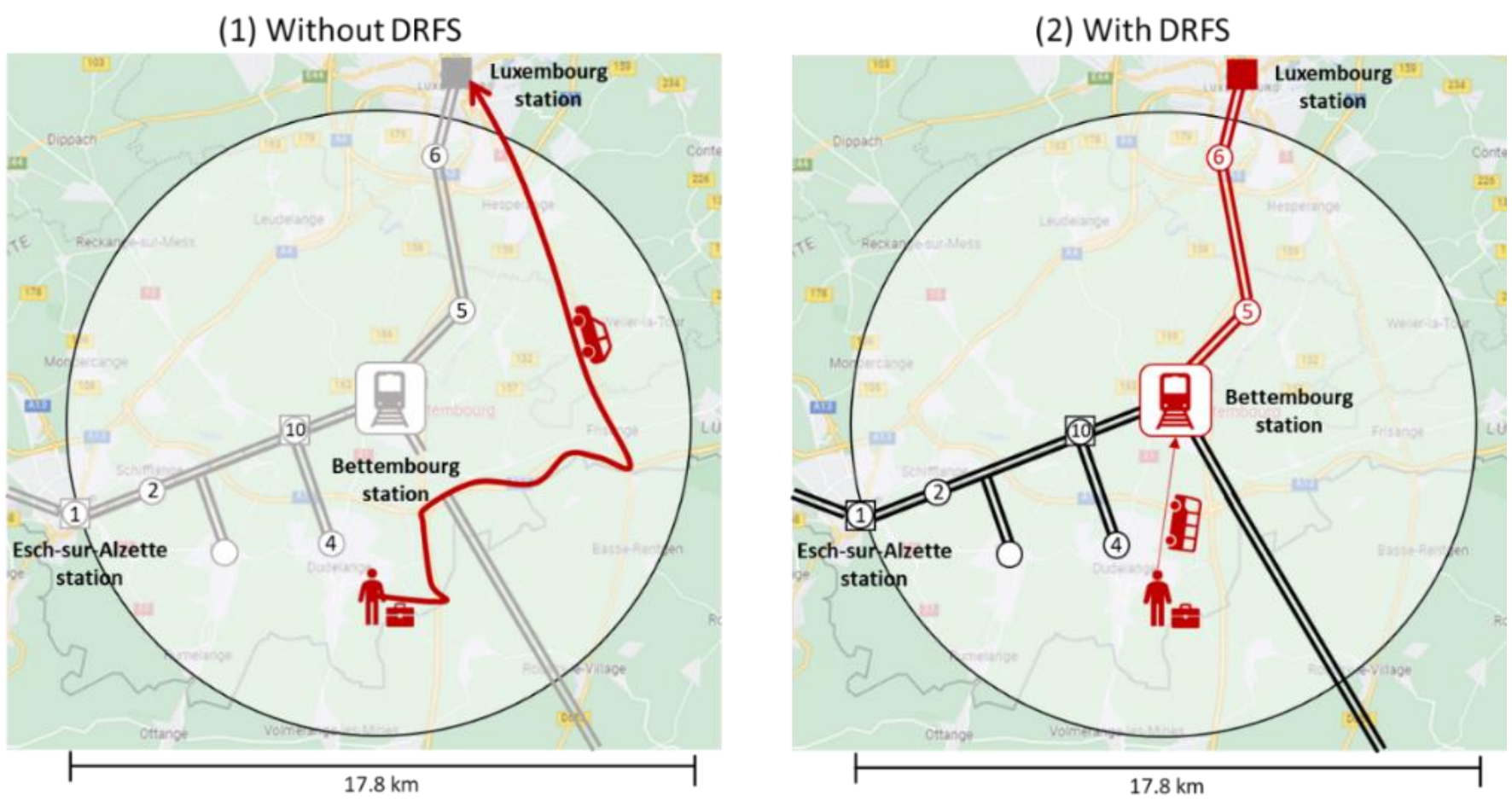

**Figure 5. An example of people replacing their car trip to Luxembourg station with the first-mile service and train.**

Two types of electric vehicles (EVs) and one type of gasoline vehicle (GV) are considered for the vehicle configurations (see Table 2). In reality, trains depart from Bettembourg station between 4:00 and 24:00, with varying frequencies throughout the day, as illustrated in Figure 6. It is important to note that Figure 6 summarizes the departure times for all train lines, irrespective of their direction.

We assume 200 requests are made throughout the day with the desired pick-up or drop-off bus stops and the departure/arrival time of the train that users are willing to take. The desired departure time of the train for each user is randomly assigned based on the train departure times from Bettembourg station between 6:00 AM and midnight. Consequently, the temporal distribution of requests throughout the day reflects the frequency of trains at Bettembourg station. The potential vehicle stops are uniformly distributed in the service area at a 1km distance apart. Walking distance from their origin and destination to the pickup/drop-off is ignored in this study as it will be around 10 minutes maximum, and it won't happen a lot. Each request contains one to four passengers. This service assumes the pre-booking system so that the schedule won't be updated dynamically during the day. The maximum user waiting time at the train station (transfer buffer time) is set as 10 minutes. The impact of this time on user inconvenience and operational costs is detailed in Appendix A3. The maximum in-vehicle time is also set as 1.5 times the direct travel time between customers' pickup and drop-off locations to account for a maximum acceptable detour time for the passengers. Service time at each pickup location is assumed to be 30 seconds. The entire fleet is fully charged at the beginning of the service. Note that we use the parameters that were tuned in Section 4.2. A transferability analysis of these parameters on the Bettembourg case study shows that they are highly transferable.

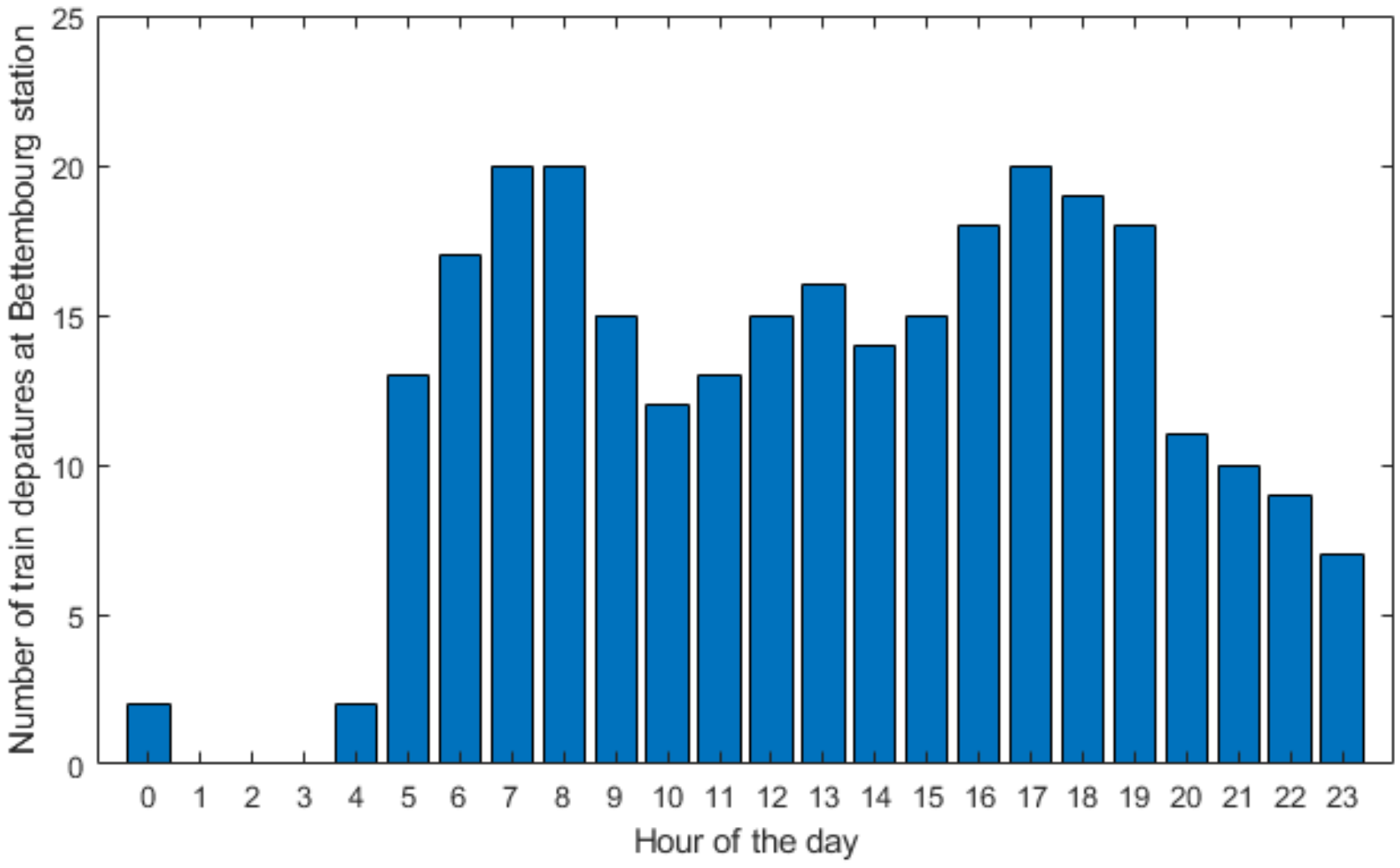

**Figure 6. The number of train departures from the Bettembourg station, including all 3 train lines, on a normal weekday.**

4.4.2. Scenario description

In this study, an opportunity charging station is assumed to be located at the Bettembourg station. Since the train station is always included as either the origin or destination, it serves as a logical and efficient location for the charging station. Additionally, the pickup and drop-off locations for each request are uniformly and randomly distributed across the area, further justifying the station's central role in the operation. To investigate the impact of charging infrastructure configuration, eight scenarios are investigated, varying by vehicle type, charger power, and $CO_2$ reduction target, as summarized in Table 6. Note that we exclusively consider high-power chargers for opportunity charging because they are most relevant for the fleet used for the DRT service, which usually consists of minibuses with 10-30 seats. Slower chargers might not efficiently recharge the fleet, resulting in a high turnaround time and the need for more vehicles to meet customer demand. For each scenario, seven cases are created, with the number of chargers ranging from 0 to 6. Note that we allow a total maximum of 30-minute waiting times for the charging operations of vehicles during the planning horizon. Our preliminary analysis shows that increasing the maximum waiting times to 1 hour leads to similar optimal solutions. This flexibility allows operators to consider the trade-offs between the investment costs and waiting costs (utilisation) of charging infrastructure. Additionally, a reference case using only GVs is included for comparison.

- 5 demand scenarios have been created, and all cases are tested with each demand scenario to compute the expected mix fleet routing costs
- For all scenarios, the configuration of GVs remains the same as presented in Table 2
- Table 7 shows the results of GV-only cases. Maximum targeted CO2 emissions are estimated for each demand scenario based on using a fleet of GVs.

**Table 6. The scenario specification.**

| Type of vehicle | Charging power | $CO_2$ reduction target | # of chargers |
|---|---|---|---|
| Type 1 | 125 kW | 90% | 0-6 |
| (24-seat) | | 50% | 0-6 |
| | Super-fast | 90% | 0-6 |
| | | 50% | 0-6 |
| Type 2 | 125 kW | 90% | 0-6 |
| (10-seat) | | 50% | 0-6 |
| | Super-fast | 90% | 0-6 |
| | | 50% | 0-6 |

Remark: The detailed characteristics of vehicles and chargers are described in Table 2.

**Table 7. Total costs, fleet size, and $CO_2$ emission using a fleet of GVs**

| Demand scenario | Total cost | GV fleet size | $CO_2$ emission |
|---|---|---|---|
| 1 | 122.17 | 5 | 268.532 |
| 2 | 123.45 | 5 | 276.818 |
| 3 | 127.94 | 5 | 305.975 |
| 4 | 110.69 | 4 | 298.972 |
| 5 | 120.53 | 5 | 257.871 |

4.4.3. Results

Table 8 summarises the optimal configuration of the charging facilities and fleet size. In general, for the cases with 24-seat vehicles, the number of chargers has minimal influence on fleet size and operational costs. This aligns with expectations, as these vehicles have larger battery capacities compared to 10-seat vehicles. The effect is particularly evident in scenarios with a 50% $CO_2$ reduction target, where half of the demand can still be met by conventional gasoline vehicles (GV), reducing the reliance on EV charging infrastructure.

It should be noted that for the cases involving 10-seat vehicles with a 90% $CO_2$ reduction target, the fleet size varied considerably among the five demand scenarios. Notably, even when three rapid chargers were available, 36 EVs were required—nearly identical to the 37 EVs needed in the scenario without any chargers. This similarity arises because, for two out of the five demand scenarios, the required EV fleet size did not start to decline significantly until a fourth rapid charger was added. Although the optimal EV fleet size for the remaining three scenarios decreased by more than 10 vehicles when additional chargers were introduced, ensuring that 90% of the demand is met across scenarios necessitates satisfying four out of the five demand cases. Consequently, the EV fleet size remained at 36 vehicles in the optimal configuration.

Figure 7 illustrates the mean operational cost without the price of chargers (top left), mean total cost with the price of chargers (bottom left), mean EV fleet size (top middle), and mean GV fleet size (bottom middle) and the mean charging time (top right) across five demand scenarios for a case using a 24-seat vehicle with a 90% $CO_2$ reduction target. Each sub-figure presents two cases with different chargers—super-fast chargers (220 kWh) and rapid chargers (125 kWh)—distinguished by different colours (blue and orange) Figure 8 illustrates the same set of results for the case with a 50% $CO_2$ reduction target.

According to Figure 7, even in the case with a 90% CO2 reduction target, charging activities are relatively limited, with a maximum observed (i.e. less than 30 minutes mean charging time among the mean fleet size of 8). In addition, the operational cost does not decrease significantly as the number of chargers increases. Since the installation cost of chargers rises with the number of chargers, the total cost—including both operational and installation costs—ultimately increases as more chargers are added, particularly in the case of super-fast chargers. For rapid chargers, however, the total cost is minimized when a single charger is installed. This is primarily due to the lower installation cost of rapid chargers, making them a more cost-effective option in this scenario.

Unlike the case with 24-seat vehicles, the EV fleet size decreases as the number of chargers increases, leading to a reduction in operational costs as well. Figure 9 and Figure 10 illustrate the same analysis with Figure 7 but for the case of a 10-seat vehicle, considering two $CO_2$ reduction targets: 90% and 50%.

For the 90% $CO_2$ reduction target, the mean charging time initially increases with the number of chargers, peaking at around 5 chargers for rapid charging and 4 chargers for super-fast charging (see Figure 9). This suggests that adding an extra charger significantly impacts charging efficiency up to these thresholds. A similar trend is observed for the mean EV fleet size and operational costs. For the 50% $CO_2$ reduction target, the operational cost steadily decreases until the number of chargers reaches 2 for both rapid and super-fast chargers (see Figure 10). Beyond this point, further increasing the number of chargers does not significantly reduce operational costs, which remain within a similar range. As seen in the top-right figure, the increase in mean charging time slows down after 2 chargers, along with the rate of EV fleet size reduction.

**Table 8. The optimal configurations for charging facilities and mixed fleet.**

| Type of vehicle | Charging power | CO2 reduction target | Mean operational costs | Mean total cost | Optimal number of chargers | Fleet configurations | | cpu (sec) |
|---|---|---|---|---|---|---|---|---|
| | | | | | | EV | GV | |
| Type 1 (24-seat) | Rapid | 90% | 429.93 | 443.63 | 1 | 7 | 1 | 300 |
| | | 50% | 306.12 | 306.12 | 0 | 4 | 3 | 161 |
| | Super-fast | 90% | 443.63 | 443.63 | 0 | 8 | 1 | 109 |
| | | 50% | 306.12 | 306.12 | 0 | 4 | 3 | 161 |
| Type 2 (10-seat) | Rapid | 90% | 943.72 | 976.60 | 3 | 36 | 1 | 1751 |
| | | 50% | 543.51 | 566.80 | 2 | 11 | 4 | 2192 |
| | Super-fast | 90% | 991.09 | 991.09 | 0 | 37 | 1 | 211 |
| | | 50% | 604.57 | 604.57 | 0 | 20 | 3 | 482 |

For both cases, while the introduction of super-fast chargers reduces operational costs, the savings do not outweigh the high installation costs of these chargers. Consequently, the scenario without super-fast chargers is found to be optimal, as shown in Table . This is the same for the case with 24-seat vehicles. It is important to note that the installation costs of super-fast chargers vary greatly depending on factors such as location (city or country) and usage (private, shared, or public). In some cases, introducing super-fast chargers could still be cost-effective when considering both operational and charging facility costs. In contrast, rapid chargers emerge as a more cost-effective solution despite their slightly higher operational costs.

Overall, rapid chargers provide the best balance between cost and efficiency for 10-seat EVs, while super-fast chargers are not cost-effective despite their operational advantages in the current settings. Table 9 summarizes the operational cost savings for scenarios with 1 to 6 super-fast chargers, compared to the no-charger scenario, under $CO_2$ reduction targets of 50% and 90%. According to Table 9, if the daily equivalent cost of a super-charger is 85.14% of the current setting, the optimal number of super-chargers would be 2 when $CO_2$ reduction target is 50%. However, if the cost is reduced to 60% of the current setting, the optimal number increases to 4. When $CO_2$ reduction target is 90%, if the daily equivalent cost of super-charge is 74% of current cost settings, the optimal number of chargers would have been 5. For 24-seat vehicles, the installation of additional chargers has a negligible impact on fleet size due to their higher battery capacity. The results demonstrate the developed model and algorithm could support the operator to optimise their fleet composition and charging infrastructure during the transition period of their fleet electrification. Other scenarios (demand, types of fleet and charging technologies) can be further studied to adapt the operator's needs in practice, which remains future extensions of this study. Intermediate targets (e.g., 70%) would produce more nuanced and transitional results. This topic is left for future exploration of this work.

**Table 9. Operational cost savings for various cases with different numbers of super-fast chargers compared to having no charging station**

| Number of super-fast chargers | | 1 | 2 | 3 | 4 | 5 | 6 |
|---|---|---|---|---|---|---|---|
| Amount of operational cost savings compared to not having a charging station | 50% $CO_2$ reduction | 32.396 | 90.621 | 68.285 | 86.514 | 80.283 | 82.854 |
| | 90 %$CO_2$ reduction | 0 | 24.812 | 62.766 | 95.002 | 121.835 | 129.443 |
| Cost of rapid chargers in this study's setting | | 13.799 | 23.488 | 33.177 | 42.866 | 52.555 | 62.244 |
| Cost of super-fast chargers in this study's setting | | 87.258 | 106.436 | 125.614 | 144.792 | 163.970 | 183.148 |

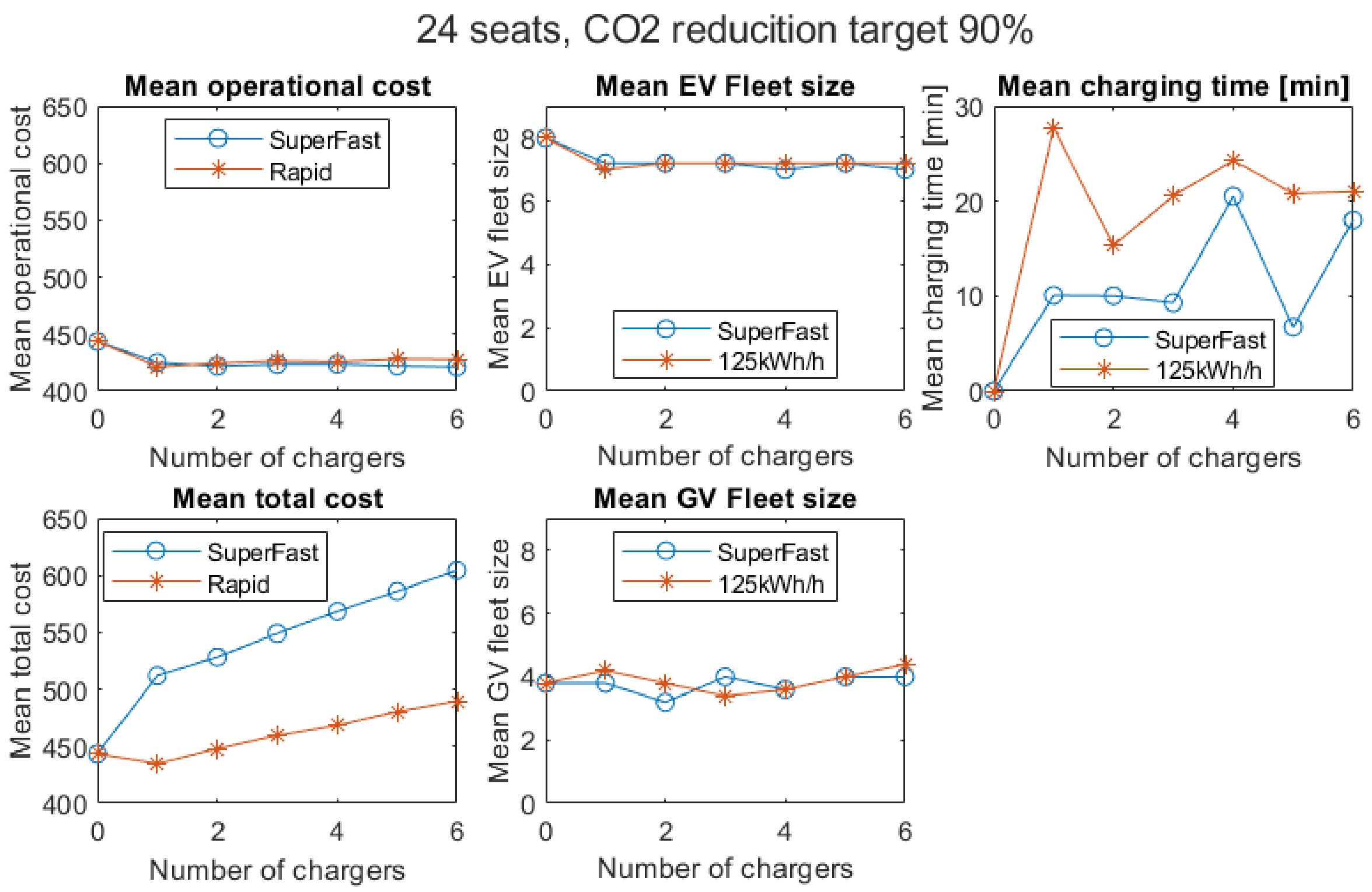

**Figure 7. The key statistics for cases with 24-seats vehicle with 90% $CO_2$ reduction target with super-fast chargers (blue) and rapid chargers (orange).**

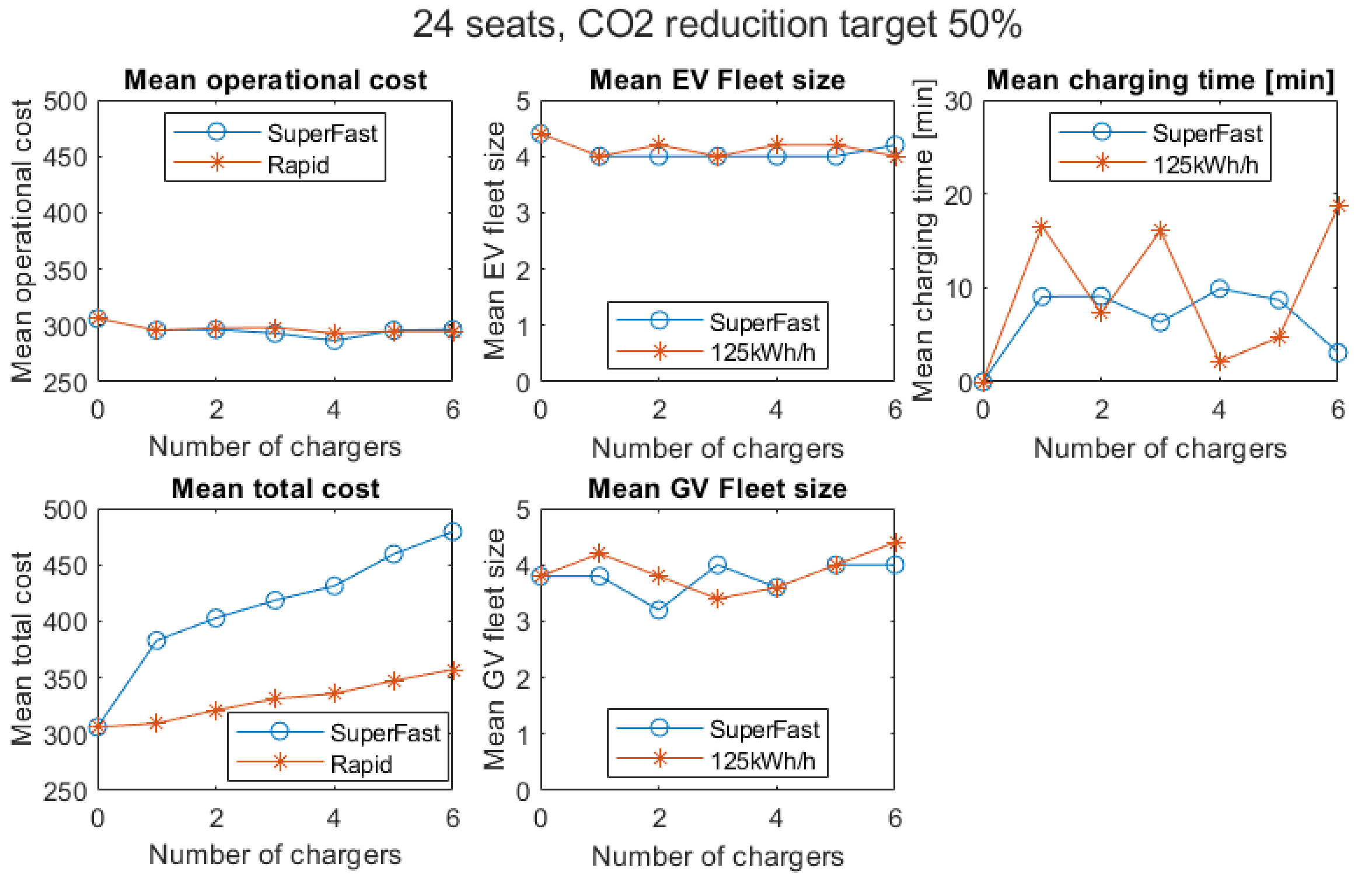

**Figure 8. The key statistics for cases with 24-seats vehicles with a 50% $CO_2$ reduction target with super-fast chargers (blue) and rapid chargers (orange).**

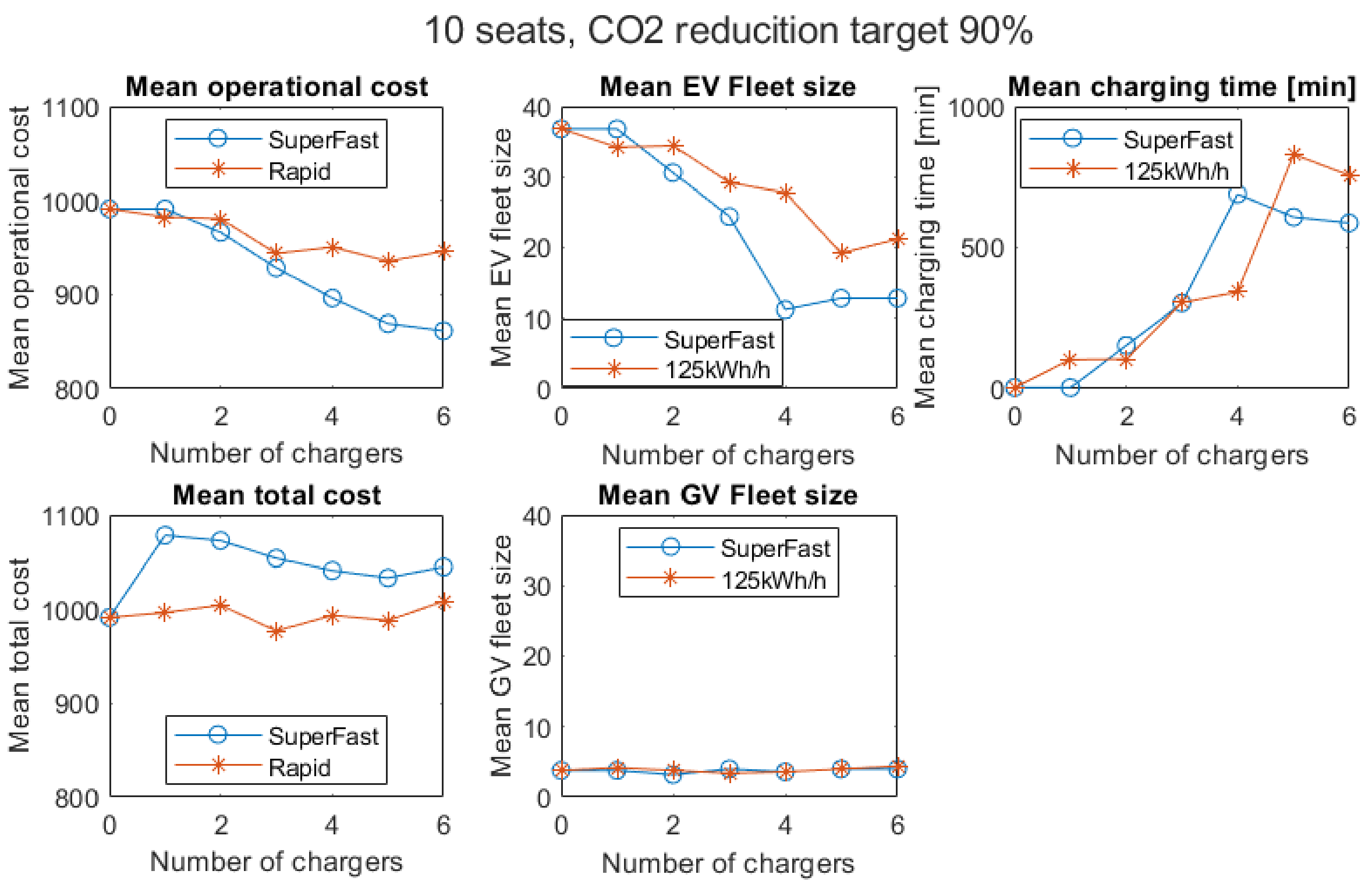

**Figure 9. The key statistics for cases with 10-seats vehicles with a 90% $CO_2$ reduction target with super-fast chargers (blue) and rapid chargers (orange).**

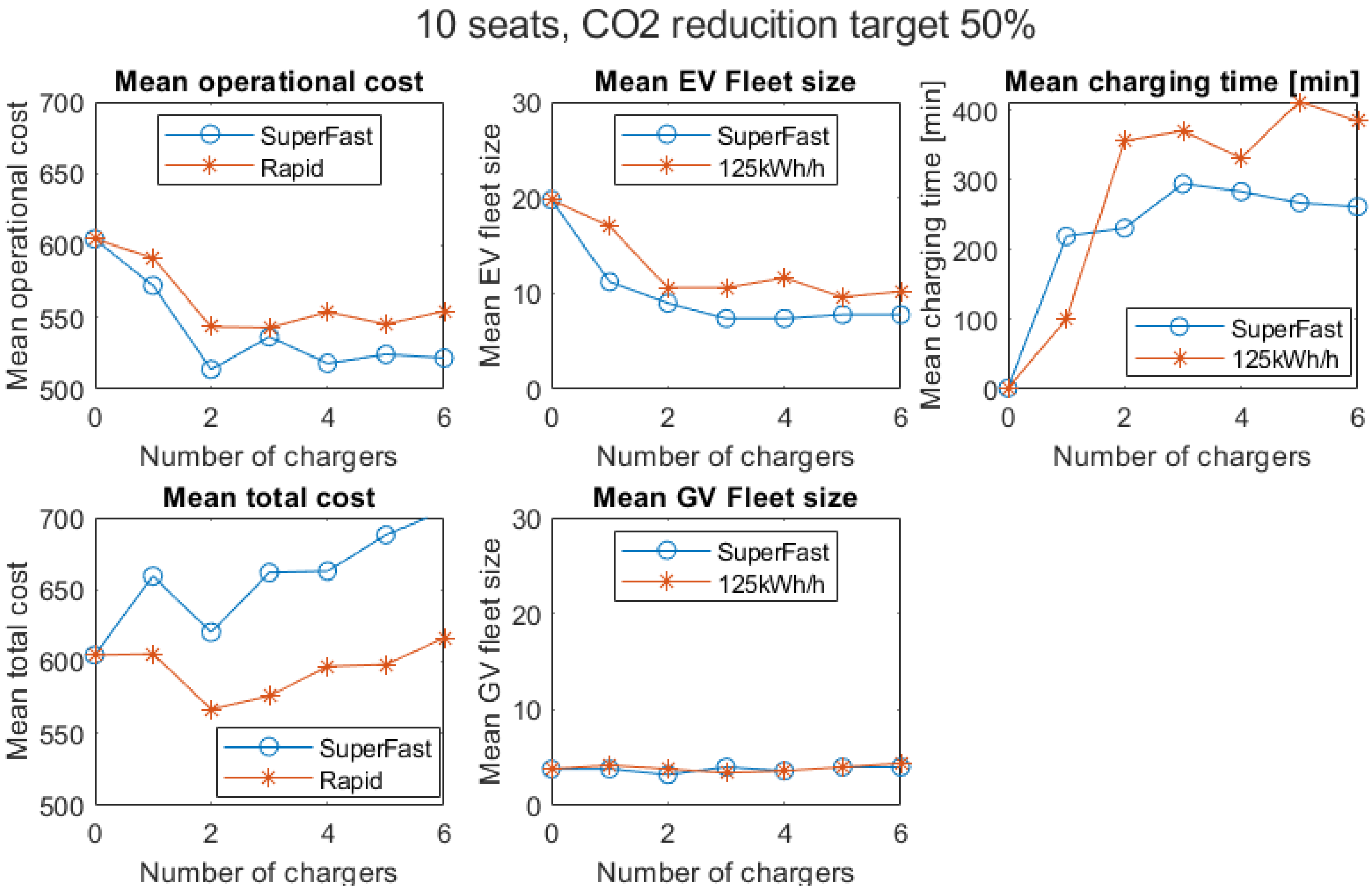

**Figure 10. The key statistics for cases with 10-seats vehicles with a 50% $CO_2$ reduction target with super-fast chargers (blue) and rapid chargers (orange).**

## 5. Conclusions and discussion

Jointly considering the trade-off of the $CO_2$ reduction target and the fleet electrification investment cost is an important decision issue for many mobility service operators as the investment costs for fleet electrification are very expensive. Existing studies assume a fully electrified fleet and jointly optimise fleet composition and charging infrastructure planning. However, the transition to fully-electrified mobility systems needs to develop a flexible modelling approach to support this process. Few studies address this research gap to optimise the electrified mobility system configuration with the presence of conventional gasoline-powered vehicles and electric vehicles with customisable $CO_2$ emission reduction targets.

This study enables this configurable $CO_2$ reduction target as a lever for the joint fleet size and charging infrastructure planning of demand responsive feeder service using a mix fleet. A bi-level optimisation model is developed to minimise the expected system costs, incorporating the charging infrastructure and mix fleet composition under stochastic customer demand and customised $CO_2$ reduction targets. The bi-level problem needs to iteratively solve a series of the $CO_2$-constrained mix fleet dial-a-ride problem with different candidate charging infrastructure configuration to obtain an optimised system configuration. The optimization problem is particularly challenging due to the need to simultaneously address multiple, interdependent objectives and constraints: 1) Minimizing the mixed fleet size and vehicle routing costs while ensuring efficient operations; 2) Strictly adhering to the CO2 emission targets, which limits the allowable emissions from gasoline vehicles; 3) Managing charging schedules effectively under limited charging station capacity, which can create bottlenecks affecting overall service performance; and 4) Balancing complex interdependencies between operational constraints (such as energy consumption, charging availability, and CO2 emissions) and customer-related requirements (such as meeting time windows and minimizing maximum ride time detours). These interconnected factors make the optimization process computationally demanding and difficult to solve efficiently. We develop a deterministic annealing-based metaheuristic with tailored penalty functions and local search operator design to obtain good solutions efficiently (around 3 minutes or less) with up to 500 requests. The proposed algorithm is validated by comparing it with the benchmarks obtained by the commercial mixed integer programming solver on a series of test instances and scenarios.

We further apply the model and algorithm for a real-world case study in Bettembourg, Luxembourg. This case study considers the electrification of a feeder service connecting the main train station in Bettembourg. Different scenarios are considered with respect to different types of EVs (with one type of GVs), charging power (rapid (125kW) and super-fast (220 kW)), and target $CO_2$ emission targets. The results demonstrate that the proposed model allows an optimised system configuration under stochastic demand.

In our analysis (see Table 2) we assumed zero $CO_2$ emissions per kilometre for EVs and consider operational (tank-to-wheel) $CO_2$ for gasoline buses (GBs). While this simplifying assumption isolates the operational benefit of electrification, it ignores the carbon intensity of grid electricity. Climatiq (2023) reports a total supplier mix of 172.5 g $CO_2$e/kWh for Luxembourg, which will decrease with ongoing investment in renewables. With consumption at 0.938 kWh/km, the BEB emits ≈162 g$CO_2$e/km.

However, our baseline assumption of 176 g$CO_2$/km for the 24 seat diesel bus is highly conservative; realistic fuel consumption of such a bus in urban service would realistically give emissions of 400–600 g$CO_2$/km. For a fair comparison we should include the full carbon cost of diesel. The JEC Well-to-Wheels Report (Prussi et al., 2020) indicates upstream Well-to-Tank (WTT) emissions for conventional diesel add approximately 25% to the operational Tank-to-Wheel (TTW) emissions. Applying this +25% factor to our diesel $CO_2$ emission figures would be required for a direct comparison with the full carbon impact of electricity. Ultimately it is beyond the scope of this paper to do a comprehensive life-cycle comparison. Recent studies for the EU found that the full life-cycle emissions of electric vehicles EU are 73% lower for cars (Negri and Bieker 2025) and 63% lower for buses (O'Connell et al 2023) than gasoline vehicles; demonstrating the substantial advantage of electrification. A natural extension of this work would be to consider varying such grid emission factors reflecting future decarbonisation pathways.

Future research directions include scenario-based analysis and system evaluation with respect to key system parameters (number of requests, detour factor, multiple depots, etc.). Another interesting extension is considering the multi-period strategic planning problem to achieve target $CO_2$ emission goals over different years. Moreover, for more complicated system configurations involving a higher number of choice options (types of vehicles, chargers, locations, etc.), the proposed metaheuristic can be further improved with additional tailored local search operators. Using a surrogate-based optimisation approach is another efficient approach to configure large systems. The influence of varying user request patterns on parameter stability can be further investigated. A comprehensive life-cycle cost assessment, including investment recovery period and operational expenses, is valuable for evaluating the economic viability and providing practical insights for real-world implementation. Furthermore, our approach can be further extended to include the joint planning with the power system (e.g., time-varying grid capacity integration) and the time-of-use energy prices for smart charging management.

**Acknowledgements**

The work was supported by the Luxembourg National Research Fund (C20/SC/14703944).

## Appendix A1. Summary of the literature review.

| Ref | Fixed-route and timetable | Demand | Fleet | Charging policy | CS capacity | Modelling & solution algorithm | Fleet size | CS location and equipment decision | Routing decision | $CO_2$ constr. |
|---|---|---|---|---|---|---|---|---|---|---|
| Chen et al. (2016) | SAEV | Stochastic | EV | Fully | | Simulation | ✓ | CS location & # of chargers | ✓ | |
| Rogge et al. (2018) | ✓ | Fixed | EB | Fully | ✓ | MILP & GA | ✓ | CS location & # of chargers | | |
| An (2020) | ✓ | Stochastic | EB | Fully | ✓ | Stochastic IP & Lagrangian relaxation | ✓ | CS location only | | |
| Zhang et al. (2020) | SAEV | Stochastic | EV | Fully | ✓ | MIP & heuristic | ✓ | CS location & # of chargers | | |
| Guschinsky et al. (2021) | ✓ | Fixed | EB + GB | Fully | ✓ | MILP & PSO combined heuristic | ✓ | CS location & # of chargers | | |
| Hiermann et al., 2016 | Mix VRPTW | Fixed | EV | Fully | | MILP & ALNS | ✓ | | ✓ | |
| Lee et al. (2021) | ✓ | Fixed | EB | Fully | ✓ | MIP & two-stage algorithm | ✓ | # of chargers at the depot only | | |
| Ma and Xie (2021) | Ride-hailing | Stochastic | EV | Partial | ✓ | Bilevel MILP & Surrogate optimisation | | CS location & # of chargers | ✓ | |
| Yildrim and Yildiz (2021) | ✓ | Fixed | EB | Fully | ✓ | MILP & CG | ✓ | CS location & # of chargers | | |
| Alam and Guo (2022) | Ride-hailing | Fixed | EV | Fully | | MILP | | CS location & # of chargers | ✓ | |
| Hu et al. (2022) | ✓ | Fixed | EB | Partial | ✓ | Robust MILP | # of EB batteries | CS location & # of chargers | | |
| Li et al. (2022) | Taxi | Stochastic | EV | Partial | ✓ | Simulation & PSO | | CS location & # of chargers | ✓ | |
| Li et al. (2022) | DRT | Fixed | GV | NA | NA | MILP | NA | NA | ✓ | ✓ |
| Paudel and Das (2022) | SAEV | Stochastic | EV | Fully | | Robust MILP | ✓ | CS location & # of chargers | | |
| Soltanpour et al. (2023) | ✓ | Fixed | EB + GB | Fully | ✓ | MILP & SA | ✓ | CS location & # of chargers | | |
| Yang et al. (2023) | SAEV | Stochastic | EV | Fully | ✓ | Simulation & BO | ✓ | CS location & # of chargers | ✓ | |
| Alamatsaz et al. (2024) | ✓ | Fixed | EB | Fully | ✓ | MILP and B&P | ✓ | CS location | | |
| Nath et al. (2024) | ✓ | Fixed | EB | Partial | ✓ | MILP & ILS | ✓ | CS location & # of chargers | | |
| Sen et al. (2024) | ✓ | Fixed | EB + GB | Fully | ✓ | MIP & heuristic | ✓ | | | |
| Shi et al. (2025) | Logistics | Fixed | EV + GB | Fully | | MILP & hybrid GA | | | ✓ | |
| Zhang et al. (2025) | SAEV | Stochastic | EV | Partial | ✓ | MILP | ✓ | CS location & # of chargers | | |
| **This study** | **Dial-a-Ride feeder service** | **Stochastic** | **EB + GB** | **Partial** | ✓ | **Bilevel MILP & DA** | ✓ | **CS location & # of chargers** | ✓ | ✓ |

Remark: **ALNS**: Adaptive Large Neighbourhood Search; **BO**: Bayesian Optimisation; **CG**: column generation; **DA**: deterministic annealing metaheuristic; **SA**: Simulated annealing; **ILS**: iterated local search; **B&P**: Branch-and-Price; **EB**: electric buses; **Fully**: charge to battery capacity or a predefined maximum level (e.g., 80% of battery capacity); **GA**: genetic algorithm; **GB**: gasoline buses; **IP**: integer program; **Mix VRPTW**: Vehicle Routing Problem with Time windows using a fleet of EVs of different types; **PSO**: Particle Swarm Optimization; **SAEV**: shared autonomous electric vehicle. **NA**: not applicable.

**Appendix A2 Transferability of the algorithmic parameters**

In this appendix, we test the transferability of the algorithmic parameters on the sets of small test instances with 10,20,…,50 requests (Table 4) and that of the Bettembourg case study. For the small test instances in Table 4, three randomly generated datasets are used (dataset 1 to dataset 3), each containing 5 test instances with 10, 20, 30, 40, and 50 requests. One dataset is used to tune the parameters of the DA metaheuristic and then we apply the tuned parameters to solve the other 10 test instances (datasets 2 and 3). The tuned parameters are referred to as Parameters T1. For the Bettembourg case study, one test instance is used to tune the parameters (referred to as Parameters T2) and apply them to the other out-of-sample 5 demand scenarios (test instances), each with 200 requests and one out-of-sample scenario is newly generated.

We compare the performance of the target-tuned parameters (T1 and T2) with that of the baseline parameters (tuned in Section 4.2 using 5 randomly generated datasets with 20,40,60,80, and 100 requests and an initial SoC of 50% for vehicles). The baseline parameter setting is referred to as T0. The parameter turning follows the same procedure described in Section 4 with an initial SoC of 50% for vehicles and 50% $CO_2$ emission reduction target. The computational time budget parameter is fixed as $Iter_{max} = 10^5$ iterations for all three parameter settings.

The target-tuned parameters, T1 and T2, are shown in Table A2. We can observe that there are significant differences between the three sets of tuned parameters. We used T1 to solve the test instances in datasets 2 and 3 with 90% and 50% $CO_2$ reduction and all GVs. The results are reported in Table A3. For the -90% and -50% $CO_2$ reduction cases, we observe that using parameters T1 achieves better objective function values: -1.37% and -0.3%, respectively, compared to using parameters T0. The p-values from the Wilcoxon signed-rank test show that there are no significant differences in performance (at a 5% level of statistical significance) when using T0 or T1. For the Bettembourg case study, using target-tuned parameters T1 achieves better performance, with an average of -0.41%. At a 5% level of statistical significance, the Wilcoxon signed-rank test shows that there are no significant differences in DA metaheuristic performance using parameters T0 and T2.

**Table A2. The baseline parameters and the tuned parameters using the target datasets.**

| Tuned parameters | Parameter values |
|---|---|
| T0 (baseline) | 100, 200, 20, 0.3, 250, 150, 300, 400, 0.9 |
| T1 (c10-c50) | 60, 300, 20, 0.1, 200, 200, 100, 400, 0.6 |
| T2 (Bettembourg) | 40, 300, 40, 0.1, 100, 200, 200, 100, 1.5 |

Note: The parameter values are referred to $\rho_1, \rho_2, \rho_3, \rho_4, \bar{n}_{remove}, \bar{n}_{stagnant}, n_{imp}, T_{red}, T_{max}$. Recall that $\rho_1$ is the penalty for each unserved user; $\rho_2$ is the penalty for each used vehicle; $\rho_3$ is the penalty for excess CO2 emission compared to the targeted threshold and $\rho_4$ is the weighted coefficient for each kWh charged to satisfy the energy constraint.

**Table A3. Performance comparison using the baseline parameters T0 and target-tuned parameters T1.**

| | | CO2 reduction target | | | | GV | |
|---|---|---|---|---|---|---|---|
| | | -90% | | -50% | | | |
| | Test instance | Para. T0 | Para. T1 | Para. T0 | Para. T1 | Para. T0 | Para. T1 |
| | | Ave. Obj. | | Ave. Obj. | | Ave. Obj. | |
| Dataset 2 | c10 | 50.66 | 50.66 | 50.66 | 50.66 | 19.54 | 19.54 |
| | c20 | 92.35 | 92.35 | 65.43 | 65.27 | 37.97 | 37.97 |
| | c30 | 148.50 | 148.50 | 100.40 | 100.05 | 58.20 | 58.20 |
| | c40 | 161.48 | 160.67 | 115.27 | 114.60 | 59.54 | 59.50 |
| | c50 | 212.87 | 205.11 | 131.87 | 130.88 | 63.00 | 62.99 |
| Dataset 3 | c10 | 75.35 | 75.35 | 56.33 | 56.33 | 35.83 | 35.83 |
| | c20 | 98.89 | 98.89 | 69.12 | 69.00 | 38.68 | 38.68 |
| | c30 | 135.44 | 130.74 | 95.08 | 94.69 | 56.79 | 56.79 |

| | | | | | | |
|---|---|---|---|---|---|---|
| c40 | 175.51 | 174.60 | 119.45 | 118.54 | 60.38 | 60.37 |
| c50 | 206.95 | 202.49 | 125.48 | 126.32 | 46.78 | 46.78 |
| Average | 135.80 | 133.94 | 92.91 | 92.63 | 47.67 | 47.66 |
| Mean diff. (T1-T0) | -1.863 (-1.37%) | | -0.275 (-0.3%) | | 0.0064 (+0.013%) | |
| Wilcoxon p-value | 0.0625 | | 0.0742 | | 0.125 | |

Remark: The average objective function values are based on 5 runs

**Table A4. Performance comparison using the baseline parameters T0 and target-tuned parameters T2.**

| | | **Para. T0** | **Para. T1** |
|---|---|---|---|
| $CO_2$ target | Scenario | Ave. Obj. | Ave. Obj. |
| 10% | 2 | 410.51 | 407.74 |
| | 3 | 466.33 | 464.92 |
| | 4 | 428.07 | 429.77 |
| | 5 | 398.59 | 401.22 |
| | 6 | 464.38 | 460.60 |
| 50% | 2 | 292.71 | 291.54 |
| | 3 | 326.02 | 325.60 |
| | 4 | 291.11 | 294.38 |
| | 5 | 276.47 | 272.45 |
| | 6 | 310.40 | 312.77 |
| GV | 2 | 122.80 | 122.39 |
| | 3 | 128.17 | 128.08 |
| | 4 | 110.14 | 110.00 |
| | 5 | 117.63 | 104.49 |
| | 6 | 128.24 | 128.13 |
| Average | | 284.77 | 283.61 |
| Mean diff. (T1-T0) | | -1.1665 (-0.41%) | |
| Wilcoxon p-value | | 0.2523 | |

Remark: The average objective function values are based on 5 runs

**Appendix A3. Sensitivity analysis on user inconvenience and operational costs using different transfer buffer times at transit stations**

We further explore the impact of transfer buffer times on user inconvenience and operational costs. Using a shorter transfer buffer time improves service quality by reducing waiting time, but increases operational costs due to tighter time constraints for vehicle route planning. We examine the Bettembourg case study with the 24-seat minibuses and the charging infrastructure configurations reported in Table 8. We consider the baseline and two user-defined $CO_2$ reduction targets. For the 90% $CO_2$ reduction target, one DC fast charger is installed; for the 50% $CO_2$ reduction target, no charging station is required (Table 8). We varied the transfer buffer times (i.e., 10, 15, and 20 minutes) and compared user inconvenience and operational costs. We ran 10 test instances using the same approach described in Section 4.4.1, with 200 randomly generated requests in the study area. Metrics for user satisfaction include average waiting time and the average detour factor for passengers (i.e., the ratio of passengers' ride time to their direct travel time). To evaluate operational efficiency, we consider the average number of EVs and GVs used, total vehicle kilometers traveled (KMT), and operational costs (i.e., the objective function value of the lower-level problem). The results are shown in Table A5 and Figure A1. For each $CO_2$ reduction target, we observe that increasing the transfer buffer time from 10 to 20 minutes increases passenger waiting time from around 6 to 11 minutes. There are no significant impacts on the user experience detour factor. From the operator’s perspective, using higher transfer buffer times reduces the number of vehicles used and fleet utilization costs, reducing operational costs by 4.2%-13.3%.

**Table A5. The impact of transfer buffer times on user inconvenience and operational costs.**

| CO2 reduction target | Transfer buffer time | Mean total cost | # chargers | # used EV | # used GV | Average operational costs | Ave. VKT | Ave. waiting time | Ave. ride time detour |
|---|---|---|---|---|---|---|---|---|---|
| 90% | 10 | 443.52 | 1 | 7.2 | 1.0 | 429.82 | 1873.2 | 5.9 | 1.26 |
| | 15 | 425.30 | 1 | 6.6 | 1.0 | 411.60 | 1943.7 | 8.4 | 1.22 |
| | 20 | 410.24 | 1 | 6.4 | 1.0 | 396.54 | 1829.2 | 10.8 | 1.22 |
| 50% | 10 | 312.94 | 0 | 4.6 | 2.9 | 312.94 | 2091.8 | 5.6 | 1.31 |
| | 15 | 290.83 | 0 | 4.1 | 3.0 | 290.83 | 2049.5 | 7.9 | 1.32 |
| | 20 | 271.99 | 0 | 4.0 | 2.5 | 271.99 | 2082.4 | 10.4 | 1.33 |
| All GVs | 10 | 121.80 | 0 | 0.0 | 4.8 | 121.80 | 1942.5 | 5.6 | 1.30 |
| | 15 | 109.00 | 0 | 0.0 | 4.1 | 109.00 | 1915.9 | 8.1 | 1.30 |
| | 20 | 105.58 | 0 | 0.0 | 4.0 | 105.58 | 1935.8 | 10.8 | 1.29 |

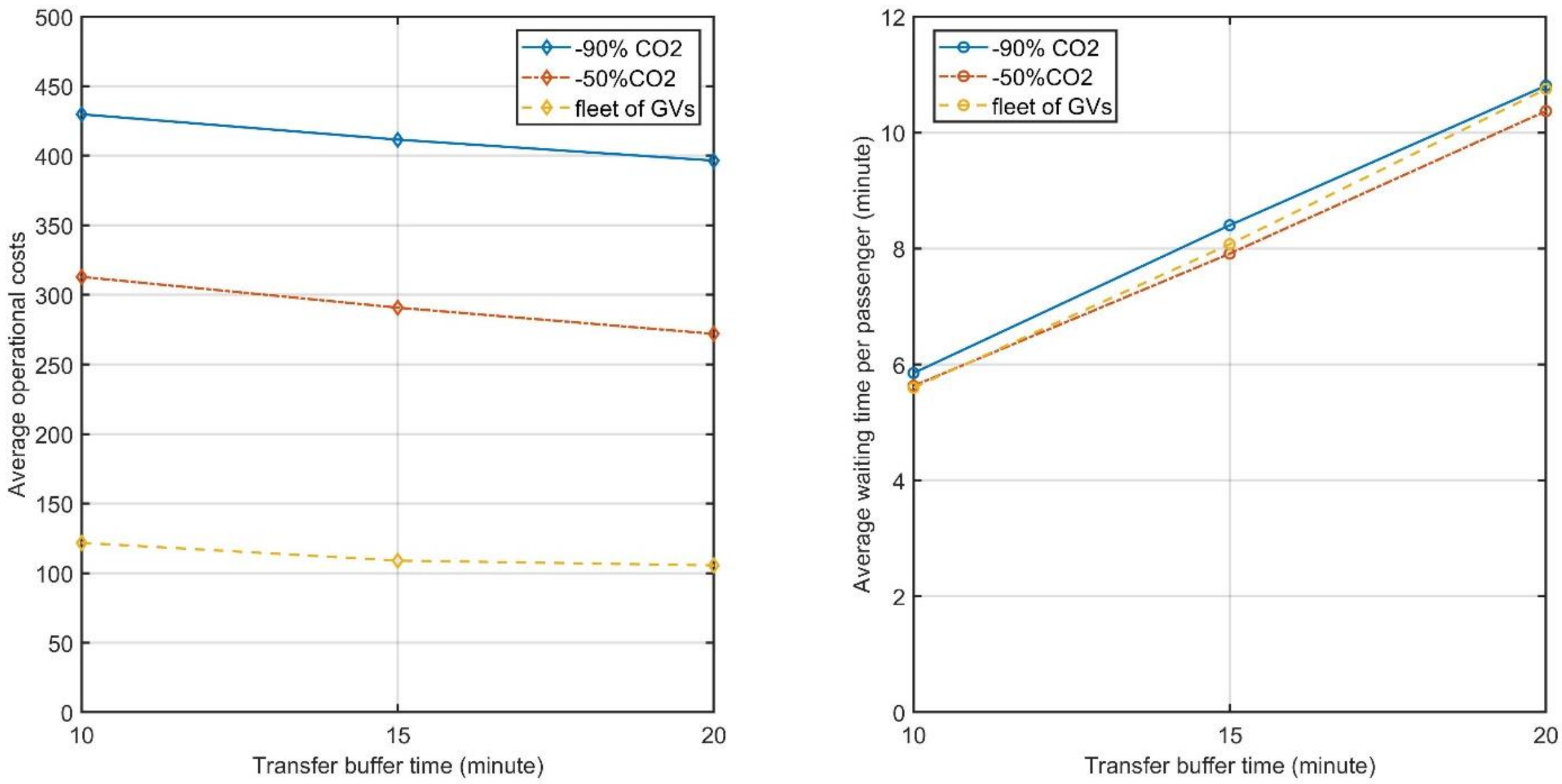

Figure A1. The impact of transfer buffer times on the operational costs (left) and the waiting times of passengers for the Bettembourg case study (24-seat minibuses) with different CO2 reduction targets.